\documentclass[a4paper,11pt,twoside,reqno,final]{amsart}

\usepackage{showkeys}

\usepackage{amsthm, color}
\usepackage{latexsym,bm}
\usepackage{bbding}

\usepackage{amsfonts}
\usepackage{amsmath}
\usepackage{mathrsfs, enumerate, amssymb}
\usepackage{tikz}
\usepackage{fancybox}
\usepackage{pstricks}
\usepackage{pst-node}
\usepackage[headinclude,DIV12]{typearea}
 \usepackage{hyperref}
\usepackage{color}
\usepackage{geometry}
\usepackage{graphicx, psfrag}

\newtheorem{theorem}{Theorem}[section]
\newtheorem{proposition}[theorem]{Proposition}

\newtheorem{lemma}[theorem]{Lemma}

\newtheorem{remark}[theorem]{Remark}
\newtheorem{example}[theorem]{Example}
\newtheorem{defi}[theorem]{Definition}

\newcommand{\me}{\mathcal{M}(E)}

\newcommand{\p}{\mathrm{P}}

\newcommand{\R}{\mathbb{R}}

\newcommand{\1}{\mathbf{1}}
\newcommand{\e}{\mathrm{e}}
\newcommand{\A}{\mathcal{A}}
\newcommand{\C}{\mathcal{C}^{\xi}_{b}(E)}
\newcommand{\ud}{\mathrm{d}}

\definecolor{wco}{rgb}{0.5,0.2,0.3}

\numberwithin{equation}{section} 

\begin{document}

\title[Law of large numbers for supercritical superprocesses with non-local branching]{Law of large numbers for supercritical superprocesses with non-local branching}
\thanks{S. P. would like to acknowledge support from a Royal Society Newton International Fellowship. The research of T.Y. is supported by NSFC (Grant No. 11501029 and 11731009).}
\author{Sandra Palau}
\author{Ting Yang }

\date{}

\maketitle
\vspace{0.2in}

\begin{abstract}

In this paper we establish a weak and a strong law of large numbers for supercritical superprocesses with general non-local
branching mechanisms.
Our results complement earlier results obtained for superprocesses
with only local branching. Several interesting examples are developed, including multitype continuous-state branching processes, multitype superdiffusions and superprocesses with discontinuous spatial motions and non-decomposable branching mechanisms.
\end{abstract}

\medskip


\medskip


\section{Introduction}\label{sec:intro}

A natural and interesting question in the theory of superprocesses is how fast the mass assigned to a compact set grows as time evolves. For superdiffusions,
Engl\"{a}nder and Turaev \cite{ET} proved a weak convergence of the ratio between the mass in a compact set and its expectation.
Later,
weak (convergence in law or in probability) and strong (almost sure convergence) laws of large numbers have been established for superdiffusions successively in \cite{EW, E,LRS,EKW} and the references therein.
For superprocesses where the spatial motion may have discontinuous paths, Chen et al. \cite{CRW} is the first paper that established the almost sure limit theorems. They showed that the principal eigenvalue of the $L^{2}$-generator associated with the mean semigroup determines asymptotic properties of the superprocesses. When the branching mechanism is purely local, the corresponding $L^{2}$-generator is a \textit{local} perturbed Schr\"{o}dinger operator (that is, the operator obtained through Feynman-Kac transform by a positive continuous additive functional).
Motivated by their work,
Wang \cite{W} and Kouritzin and Ren \cite{KR} established the
strong law of large numbers (SLLN in abbreviation)
for super-Brownian motions and super-$\alpha$-stable processes,
where the branching mechanisms are quadratic and spatially independent.
The key ingredient in their work is Fourier analysis,
which requires that the transition density of the Feynman-Kac semigroup can be represented in terms of spectral measure and the eigenfunctions of the Schr\"{o}dinger operators.
Very recently, a new approach to SLLN has been taken in \cite{CRY} and \cite{EKW}. The core of their proofs is the skeleton decomposition, that represents the (purely local branching) superprocess as an immigration process along a branching Markov process, called the skeleton. An advantage of this method is that it enables
one
to transfer results directly from the theory of branching Markov processes. However, for a general (non-local branching) superprocess, even the existence of the skeleton needs to be justified.

In the above mentioned papers, the branching mechanisms are assumed to be purely local.
Unfortunately,
there is less work on the limit theorems for non-local branching superprocesses.
In a recent paper, Kyprianou and Palau \cite{KP} established a spine decomposition for a multitype continuous-state branching process
(MCSBP in abbreviation)
and used it to study extinction properties. Concurrently to their work, a similar decomposition has been obtained by Chen et al. \cite{CRS} for a special class of multitype superdiffusions. This decomposition is further extended in \cite{RSY} to superprocess with a branching mechanism which has both local and non-local parts. Very recently, Kyprianou et al. established in \cite{KPR} the SLLN for a supercritical MCSBP. The papers mentioned above concerned only special kinds of non-local branching superprocesses. In fact, for a MCSBP (resp. a multitype diffusion), if one considers the
$E$-valued spatial motion
on an
enriched state space $E\times I$,
where $I$ is the finite or countable set of types, then the mutation in types is the jumps in the $I$-coordinates, and the associated Feynman-Kac semigroup is generated by a matrix (resp. a coupled elliptic system, cf. Example \ref{eg:KP} and Example \ref{eg:multitype superprocess} below). So, the spectral theory of matrices (resp. the potential theory for elliptic systems) can be applied.
For a general non-local branching superprocess,
the associated Schr\"{o}dinger operator takes the form
$\mathcal{J}-a+\gamma$,
where $\mathcal{J}$ is the generator of underlying spatial motion,
$a$ is a bounded function, $\gamma$ is an integral operator,  and $a$, $\gamma$ are
related
to the branching mechanism (cf. equation \eqref{mean equation} below). Since $\gamma$ can be quite general,
the methods mentioned above are not applicable and a different approach is needed.
In this paper, we characterise the Schr\"{o}dinger operator in terms of the associated bilinear form, and impose some technical assumptions (\eqref{newA1}-\eqref{newA3} below) to ensure the existence of a positive principal eigenvalue $-\lambda_{1}$ and a ground state of the Schr\"{o}dinger operator. These conditions may look strong but they hold for a large class of processes, and we illustrate this for several key examples in Section \ref{sec:examples}. Under these and a few more assumptions, we show in Theorem \ref{them:wlln} and Theorem \ref{them:slln}
that the mass of a (non-local branching) superprocess on every compact set grows exponentially at rate $-\lambda_{1}$, and the ground state determines the asymptotic distribution.
Our proof of SLLN follows two main steps, first to obtain
the
SLLN along lattice times and then
to
extend it to all times through approximation of bounded functions by resolvent functions.
Our approach to the convergence along lattice times relies on a stochastic integral representation of superprocesses (Proposition \ref{prop6.1} below). This representation enables one
to decompose the superprocess into (not necessarily orthogonal or worthy) martingale measures,
and therefore is useful in studying the structure properties of superprocesses.
We are not the first ones to use stochastic analysis to study the limit theorems of superprocesses.
A similar idea was used in \cite{KPR} for MCSBPs and in \cite{LRS} for superdiffusions on bounded domains. However, in this paper, we extend this idea much further by considering superprocesses where the spacial motion may be discontinuous and the branching mechanism is allowed to be non-local.

The remainder of this paper is organized as follows. We start Section \ref{sec:prelim} with a review on
definitions and basic properties of symmetric Borel right processes, non-local branching superprocesses,
mean semigroups and the associated bilinear forms.
In Section \ref{sec:main} we present the main results on weak and strong laws of large numbers and give concrete examples.
In Section \ref{sec:martingale and representation}, we investigate the martingale problem and establish a stochastic integral representation for superprocesses. Finally, in the last section we give the proofs of the main results.

\section{Preliminaries}\label{sec:prelim}
	
Throughout this paper, ``:=" means ``is defined to be".
Suppose that $E$ is a Luzin topological space with Borel $\sigma$-algebra $\mathcal{B}(E)$ and $m$ is a $\sigma$-finite measure on $(E,\mathcal{B}(E))$ with full support. Let $E_{\partial}:=E\cup \{\partial\}$ be the one-point compactification of $E$.
Any function $f$ on $E$ will be automatically extended to $E_{\partial}$ by
setting $f(\partial)=0$. Let $\me$ denote the space of finite Borel measures on $E$ topologized by the weak convergence
and $\me^{0}:=\me\setminus\{0\}$ where $0$ denotes the null measure. For $\mu$ a measure on $\mathcal{B}(E)$ and $f$, $g$ measurable functions, let $\langle f,\mu\rangle:=\int_{E}f(x)\mu(\ud x)$ and $(f,g):=\int_{E}f(x)g(x)m(\ud x)$ whenever the integrals make sense. Sometimes we also write $\mu(f)$ for $\langle f,\mu\rangle$.\
For a function $f$ on $E$, $\|f\|_{\infty}:=\sup_{x\in E}|f(x)|$.
If $f(x,t)$ is a function on \ $E\times [0,+\infty)$, \ we say $f$ is \textit{locally bounded} if \ $\sup_{t\in [0,T]}\sup_{x\in E}|f(x,t)|<+\infty$ \ for every \ $T\in (0,+\infty)$.\
We use \ $\mathcal{B}_{b}(E)$ \ (respectively, \ $\mathcal{B}^{+}(E)$ \ or \ $C(E)$)  to denote
the space of bounded (respectively, nonnegative or continuous) measurable functions on $(E,\mathcal{B}(E))$.
For $a,b\in \mathbb{R}$, let $a\wedge b:=\min\{a,b\}$ and $a^{-}:=\max\{-a,0\}$.

\subsection{Spatial motion}\label{sec:spatial motion}

Let \ $\xi=(\Omega,\mathcal{H},\mathcal{H}_{t},\theta_{t},\xi_{t},\Pi_{x},\zeta)$ \ be an \ $m$-symmetric Borel right process
on $E$, where \ $\{\mathcal{H}_{t}:\ t\ge 0\}$ is the associated  natural  filtration, \ $\{\theta_{t}:\ t\ge 0\}$ \ is a  time-shift operator of $\xi$ satisfying \ $\xi_{t}\circ\theta_{s}=\xi_{t+s}$ \ for $s,t\ge 0$, and \ $\zeta:=\inf\{t>0:\
\xi_{t}=\partial\}$ \ is the lifetime of $\xi$.
Denote by  $\{P_{t}:\ t\ge 0\}$
the transition semigroup of $\xi$, in other words,
 $$P_{t}f(x):=\Pi_{x}\left[f(\xi_{t})\right], \qquad \forall f\in\mathcal{B}^{+}(E).$$
 It is known that
$\{P_{t}:t\ge 0\}$ can be uniquely extended to a strongly continuous contraction semigroup on $L^{2}(E,m)$, which we also denote by $\{P_{t}:t\ge 0\}$
(cf. \cite[Lemma 1.1.14]{CF}).
Then, by the theory of Dirichlet forms, there exists a
symmetric quasi-regular Dirichlet form
$(\mathcal{E},\mathcal{F})$ on $L^{2}(E,m)$ associated with $\xi$:
 $$
 \mathcal{F}= \Big \{u\in L^{2}(E,m):\ \sup_{t>0}\ \frac{1}{t}\int_{E}\left(u(x)-P_{t}u(x)\right)u(x)m(\ud x)<+\infty \Big\},
 $$
 $$\mathcal{E}(u,v)=\lim_{t\to 0}\ \frac{1}{t}\int_{E}\left(u(x)-P_{t}u(x)\right)v(x)m(\ud x),\qquad
 \forall u,v\in\mathcal{F}.
 $$
Moreover, this process is quasi-homeomorphic to a Hunt process associated with a regular Dirichlet form on a locally compact separable metric space (cf. \cite{FOT}) and all of the results of \cite{FOT} can be applied to $\xi$ and its Dirichlet form. Henceforth,
we may and do assume that
 $\xi$ is an $m$-symmetric Hunt process on
 a locally compact separable metric space associated with a regular Dirichlet form $(\mathcal{E},\mathcal{F})$.
In addition, we assume that $\xi$ admits a transition density $p(t,x,y)$ with respect to the measure $m$, which is symmetric
in $(x,y)$ for each $t>0$.

\subsection{Non-local branching superprocesses}\label{sec:model}

In this paper,
we consider a superprocess $X:=\{X_{t}:\ t\ge 0\}$ associated to the spatial motion $\xi$ and a (non-local) branching mechanism $\psi$ given by
\begin{equation}\label{bm1}
\begin{split}
\psi(x,f):=&a(x)f(x)+b(x)f(x)^{2}-\eta(x,f)\\
&\ +\int_{\me^{0}}\left(\e^{-\nu(f)}-1 +\nu(\{x\})f(x)\right)H(x,\ud \nu),
\end{split}
\end{equation}
for $x\in E$ and \ $f\in\mathcal{B}^{+}_{b}(E)$, \ where \ $a(x)\in \mathcal{B}_{b}(E)$, \ $b(x)\in \mathcal{B}^{+}_{b}(E)$, \ $\eta(x,\ud y)$ \ is a bounded kernel on $E$ and \ $H(x,\ud \nu)$ \ is a $\sigma$-finite kernel from $E$ to $\me^{0}$ such that
$$\sup_{x\in E}\int_{\me^{0}}\left(\nu(1)\wedge \nu(1)^{2}+\nu_{x}(1)\right)H(x,\ud \nu)<+\infty.$$
Here, $\nu_{x}(\ud y)$ denotes the restriction of $\nu(\ud y)$ to $E\setminus\{x\}$.
To be specific,
$X$ is a $\me$-valued Markov process satisfying that for every $f\in \mathcal{B}^{+}_{b}(E)$ and every $\mu\in\me$,
\begin{equation*}
\p_{\mu}\left(e^{-\langle f,X_{t}\rangle}\right)=\e^{-\langle V_{t}f,\mu\rangle},\qquad\mbox{ for }t\ge 0,
\end{equation*}
where \ $V_{t}f(x):=-\log\p_{\delta_{x}}\left(\e^{-\langle f,X_{t}\rangle}\right)$ \ is the unique nonnegative locally bounded solution to the integral equation
\begin{equation*}
V_{t}f(x)=P_{t}f(x)-\Pi_{x}\left[\int_{0}^{t}\psi(\xi_{s},V_{t-s}f)\ud s\right].
\end{equation*}
Such a process is defined in \cite{Li} via its log-Laplace functional and referred to as the $(P_{t},\psi)$-superprocess.
The branching mechanisms defined in \eqref{bm1} are quite general.
For example, let
\begin{equation}\label{local part}
\phi^{L}(x,\lambda):=a(x)\lambda+b(x)\lambda^{2}+\int_{(0,+\infty)}\left(e^{-\lambda u}-1+\lambda u\right)\pi^{L}(x,\ud u),
\end{equation}
for $x\in E$ and $\lambda\ge 0$, where \ $(u\wedge u^{2})\pi^{L}(x,\ud u)$ \ is a bounded kernel from $E$ to $(0,+\infty)$, and
\begin{equation*}
\phi^{NL}(x,f):=-\eta(x,f)+\int_{\me^{0}}\left(\e^{-\nu(f)}-1\right)\pi^{NL}(x ,\ud \nu)
\end{equation*}
for $x\in E$ and \ $f\in\mathcal{B}^{+}(E)$,\  where \ $\nu(1)\pi^{NL}(x,\ud \nu)$ \ is a bounded kernel from $E$ to $\me^{0}$. Then \ $(x,f)\mapsto \phi^{L}(x,f(x))+\phi^{NL}(x,f)$ \ is a branching mechanism that can be represented in the form of
\eqref{bm1}.
A branching mechanism of this type is said to be \textit{decomposable} with local part $\phi^{L}$ and non-local part $\phi^{NL}$.
In particular, if the non-local part equals $0$, we call such a branching mechanism \textit{purely local}.
Another usual way to define superprocesses with a decomposable branching mechanism is as a scaling limit of a sequence of branching particle systems (cf. \cite{D93,D94} and \cite{Li}).

We can rewrite \eqref{bm1} into
\begin{equation}\label{1.1}
\begin{split}
\psi(x,f)&=a(x)f(x)+b(x)f(x)^{2}-\gamma (x,f)\\
&\qquad + \int_{\me^{0}}\left(e^{-\nu(f)}-1+\nu(f)\right)H(x,\ud \nu),
\end{split}
\end{equation}
where \ $\gamma(x,\ud y):=\eta(x,\ud y)+\int_{\me^{0}}\nu_{x}(\ud y)H(x,\ud \nu)$
is a bounded kernel on $E$.
We note that,
$\psi$ given by \eqref{1.1} is purely local if and only if $\gamma(x,1)=0$ for all $x\in E$.

By \cite[Theorem 5.12]{Li}, a $(P_{t},\psi)$-superprocess $X$ has a right realization in $\me$. Let denote by $\mathcal{W}^{+}_{0}$ the space of right continuous paths from $[0,+\infty)$ to $\me$ having zero as a trap. Here, we assume that $X$ is the coordinate process in $\mathcal{W}^{+}_{0}$ and  $(\mathcal{F}_{t})_{t\in[0,\infty]}$ is the filtration generated by the coordinate process, which is completed  with the class of $\p_{\mu}$-negligible measurable sets for every $\mu\in\me$. We emphasize that the branching mechanisms considered in this paper are allowed to be
non-local and non-decomposable.
In Section \ref{sec:examples} we give a concrete example of
a non-local and non-decomposable
branching mechanism (Example \ref{eg:CKS} below).

\subsection{Mean semigroups and the associated bilinear forms}\label{sec:assump}

It is known from \cite[Proposition 2.27]{Li}
that for every $\mu\in\me$ and $f\in \mathcal{B}_{b}(E)$,
\begin{equation*}
\p_{\mu}\left(\langle f,X_{t}\rangle\right)=
\langle \mathfrak{P}_{t}f,\mu\rangle,
\end{equation*}
where $\mathfrak{P}_{t}f(x)$ is the unique locally
bounded solution to the integral equation
\begin{equation}
\mathfrak{P}_{t}f(x)=P_{t}f(x)-\Pi_{x}\left[\int_{0}^{t}a(\xi_{s})
\mathfrak{P}_{t-s}f(\xi_{s})\ud s\right]
+\Pi_{x}\left[\int_{0}^{t}\gamma(\xi_{s},\mathfrak{P}_{t-s}f)\ud s\right].\label{mean equation}
\end{equation}
By the Markov property of $X$, the operator $\mathfrak{P}_{t}$ satisfies the semigroup property, i.e., $\mathfrak{P}_{t}\mathfrak{P}_{s}=\mathfrak{P}_{t+s}$ for all $t,s\ge 0$.
Moreover,
$\mathfrak{P}_{t}$ admits a transition density $\mathfrak{p}(t,x,y)$ with respect to the measure $m$.
In fact,
if $m(B)=0$ for some $B \subset E$, then by the hypothesis
$P_{t}1_B(x)=0$ for all $t\geq 0$ and $x\in E$. Therefore, $\mathfrak{P}_{t}1_B(x)=0$ is the
unique locally bounded solution
to \eqref{mean equation} for $f=1_B$. This implies that $\mathfrak{P}_{t}<<m$ and $\mathfrak{p}(t,x,y)$ exits.

We now introduce a class of nonnegative smooth measures on $E$ (cf. \cite{Chen}).
\begin{defi}
	A nonnegative measure $\mu$  on $E$ is called a smooth measure of $\xi$ if there
	is a positive continuous additive functional $A^{\mu}_{t}$ of $\xi$ such that
	$$\int_{E}f(x)\mu(\ud x)=\lim_{t\to 0}\frac{1}{t}\Pi_{m}\left[ \int_{0}^{t}f(
	\xi_{s}
	)\ud A^{\mu}_{s}\right],
	\qquad \forall f\in \mathcal{B}^{+}(E).$$
	Here, \ $\Pi_{m}(\cdot):=\int_{E}\Pi_{x}(\cdot)m(\ud x)$. \
	In this case, $\mu$ is also called the Revuz measure of $A^{\mu}_{t}$. Moreover,
	we say that a smooth measure $\mu$
	belongs to the Kato class $\mathbf{K}(\xi)$,
	if \begin{equation*}
	\lim_{t\downarrow 0}\sup_{x\in
		E}\int_{0}^{t}\int_{E}p(s,x,y)
	\mu(\ud y)\ud s
	=0.
	\end{equation*}
	A function $g$ is said to be in the class $\mathbf{K}(\xi)$ if the measure \  $g(x)m(\ud x)$ \ is in $\mathbf{K}(\xi)$.
\end{defi}
Clearly all bounded measurable functions are included in $\mathbf{K}(\xi)$.
It is known
(see, e.g., \cite[Proposition 2.1.(i)]{AM} and \cite[Theorem 3.1]{SV})
that if $\nu\in\mathbf{K}(\xi)$, then for every $\epsilon>0$ there is some constant $A_{\epsilon}>0$ such that
\begin{equation}
\int_{E}u(x)^{2}\nu(\ud x)\le \epsilon\, \mathcal{E}(u,u)+A_{\epsilon}\int_{E}u(x)^{2}m(\ud x),\qquad\forall u\in\mathcal{F}.\label{i1}
\end{equation}

First,
we assume the following condition holds.
\begin{equation}\label{newA1}\tag{\textbf{A1}}
\int_{E}\gamma(x,\ud y)m(\ud x) \mbox{ is a Kato measure of }\xi.
\end{equation}

Under condition \eqref{newA1}, it follows from \eqref{i1}, the boundedness of \ $x\mapsto\gamma(x,1)$, \ and the inequality
$$|u(x)u(y)|\le \frac{1}{2}(u(x)^{2}+u(y)^{2})$$
that for every $\epsilon>0$, there is a constant $K_{\epsilon}>0$ such that
\begin{equation}\nonumber
\int_{E}\int_{E}u(x)u(y)\gamma(x,\ud y)m(\ud x)\le \epsilon\, \mathcal{E}(u,u)+K_{\epsilon}\int_{E}u(x)^{2}m(\ud x),\qquad\forall u\in\mathcal{F}.
\end{equation}
It follows that the bilinear form $(\mathcal{Q},\mathcal{F})$ defined by
\begin{equation}\nonumber
\mathcal{Q}(u,v):=\mathcal{E}(u,v)+\int_{E}a(x)u(x)v(x)m(\ud x)-\int_{E}\int_{E}u(y)v(x)\gamma(x,\ud y)m(\ud x)
\end{equation}
for every $u,v\in\mathcal{F}$
is closed and that there are positive constants $K$ and $\beta_{0}$ such that $\mathcal{Q}_{\beta_{0}}(u,u):=\mathcal{Q}(u,u)+\beta_{0}(u,u)\ge 0$ for all $u\in\mathcal{F}$, and
$$|\mathcal{Q}(u,v)|\le K \mathcal{Q}_{\beta_{0}}(u,u)^{1/2}\mathcal{Q}_{\beta_{0}}(v,v)^{1/2},\qquad\forall u,v\in\mathcal{F}.$$
Then, from \cite{Kunita}, for the closed form $(\mathcal{Q},\mathcal{F})$ on $L^{2}(E,m)$, there corresponds a unique pair of strongly continuous, dual semigroups \ $\{T_{t}:t\ge 0\}$ and \ $\{\widehat{T}_{t}:t\ge 0\}$ on $L^{2}(E,m)$ satisfying that \ $\|T_{t}\|_{L^{2}(E,m)}\le e^{\beta_{0}t}$, \ $\|\widehat{T}_{t}\|_{L^{2}(E,m)}\le e^{\beta_{0}t}$, \ and that for all $\alpha>\beta_{0}$,
\begin{equation}\nonumber
\mathcal{Q}_{\alpha}(G_{\alpha}f,g)=\mathcal{Q}_{\alpha}(g,\widehat{G}_{\alpha}f)=(f,g),\qquad\forall f\in L^{2}(E,m),\ g\in\mathcal{F}.
\end{equation}
Here \ $G_{\alpha}f:=\int_{0}^{+\infty}e^{-\alpha t}T_{t}f \ud t$ \ and \ $\widehat{G}_{\alpha}f:=\int_{0}^{+\infty}e^{-\alpha t}\widehat{T}_{t}f \ud t$.

We make two more assumptions. Assume that
\begin{equation}\label{newA2}\tag{\textbf{A2}}
a(x),\ \gamma(x,1)\in L^{2}(E,m),
\end{equation}
and that, there
exist a constant $\lambda_{1}<0$ and strictly
positive functions \ $h,\widehat{h}\in\mathcal{F}$ \ with $h$ bounded continuous, \ $\|h\|_{L^{2}(E,m)}=1$ \ and \ $(h,\widehat{h})=1$ \ such that
\begin{equation}\label{newA3}\tag{\textbf{A3}}
\mathcal{Q}(h,v)=\lambda_{1}(h,v),\quad
\mathcal{Q}(v,\widehat{h})=
\lambda_{1}(v,\widehat{h}),\qquad\forall v\in\mathcal{F}.\nonumber
\end{equation}

It is proved in \cite{RSY} that under
\eqref{newA1}-\eqref{newA2},
for every $t>0$, $T_{t}$ is the unique bounded linear operator on $L^{2}(E,m)$ which is equal to $\mathfrak{P}_{t}$ on $L^{2}(E,m)\cap  \mathcal{B}_{b}(E)$. More precisely, for all $f\in L^{2}(E,m)\cap  \mathcal{B}_{b}(E)$, $T_{t}f=\mathfrak{P}_{t}f$ in $L^{2}(E,m)$.
On the other hand, condition \eqref{newA3} implies that $T_th=e^{-\lambda_1 t}h$ and $\widehat{T}_t\widehat{h}=e^{-\lambda_1 t}\widehat{h}$ in $L^{2}(E,m)$ for all $t\geq 0$. Therefore, conditions \eqref{newA1}-\eqref{newA3} amounts to saying that $-\lambda_{1}$ is the principal eigenvalue of the $L^{2}$-generator of the semigroup $(\mathfrak{P}_{t})_{t\ge 0}$,
and that
$h$ is the associated ground state.

Let us make a short remark on \eqref{newA3}. In the case of a purely local branching mechanism
where $\psi=\phi^{L}$ is given in
\eqref{local part},
the associated $L^{2}$-generator of $(\mathfrak{P}_{t})_{t\ge 0}$ takes the form
$\mathcal{J}-a$,
where $\mathcal{J}$ denotes the $L^{2}$-generator of underlying spatial motion.  In this case, condition \eqref{newA3} is satisfied,
for instance, by symmetric diffusions on bounded smooth domains in $\R^{d}$ as well as symmetric $\alpha$-stable processses on $\R^{d}$ (cf. Example \ref{eg:localbm} below and the references therein). In Section \ref{sec:examples} we give more examples of non-local branching superprocesses for which conditions \eqref{newA1}-\eqref{newA3} are satisfied.

\section{Main results and examples}\label{sec:main}

\subsection{Statements of the main results}\label{sec:results}

Now we are going to present the main results of this paper. The first one relates the principal eigenvalue of
$\mathfrak{P}_{t}$
and the associated ground state  with a martingale.
\medskip
\begin{proposition}\label{prop:martingale}
Suppose \eqref{newA1}-\eqref{newA3} hold.
For every $\mu\in\me$, $W^{h}_{t}(X):=e^{\lambda_{1}t}\langle h,X_{t}\rangle$ is a non-negative $\p_{\mu}$-martingale with respect to the filtration $\{\mathcal{F}_{t}:t\ge 0\}$.
\end{proposition}

We assume the following condition holds for the remainder of this paper.
\begin{equation}\label{newA4}\tag{\textbf{A4}}
\mbox{The operators }f\mapsto \psi(\cdot,f)\mbox{ and }f\mapsto \gamma(\cdot,f)-a(\cdot) f(\cdot)\mbox{  preserve }\mathcal{C}^{\xi}_{b}(E).
\end{equation}
Here $\mathcal{C}^{\xi}_{b}(E)$ denotes the set of bounded measurable functions that are finely continuous with respect to $\xi$.

Let $W^{h}_{\infty}(X)$ be the martingale limit of $W^{h}_{t}(X)$.
Our second result gives the $L^{p}$-convergence of $W^{h}_{t}(X)$
for a $p\in (1,2]$.

\begin{theorem}\label{them7.1}
	Suppose \eqref{newA1}-\eqref{newA4} hold. If there is $p\in (1,2]$ such that
	\begin{equation}\label{newA5}\tag{\textbf{A5}}
	\sup_{x\in E}\ h^{-1}(x)\int_{\me^{0}}\nu(h)^{p}H(x,\ud \nu)<+\infty,
	\end{equation}
	then, $W^{h}_{t}(X)$ converges to $W^{h}_{\infty}(X)$ in $L^{p}(\p_{\mu})$ for every $\mu\in\me$.
\end{theorem}

We define the operators
\begin{equation}
\widetilde{P}_{t}f(x)=\frac{e^{\lambda_{1}t}}{h(x)}\mathfrak{P}_{t}(fh)(x)
\label{prop2.0}
\end{equation}
and
\begin{equation}
\widetilde{p}_{t}(t,x,y)=\frac{e^{\lambda_{1}t}}{h(x)\widehat{h}(y)}\mathfrak{p}(t,x,y),\label{prop2.0 density}
\end{equation}	
for
$t\geq 0$, $x,y\in E$, $f\in \mathcal{B}^{+}_{b}(E)$,
where
$\mathfrak{p}(t,x,y)$ is the transition density of $\mathfrak{P}_t$ with respect to $m$. An intuition of the above operators is given in Section \ref{sec:interpretation}, where it is showed that they can be seen, respectively, as the transition semigroup and the transition density function (with respect to $h\widehat{h}m$) of an auxiliary process $\widetilde{\xi}$, see Proposition \ref{prop2}
and Remark \ref{remark:many to one} below.

\begin{theorem}[Weak law of large numbers]\label{them:wlln}
	Suppose \eqref{newA1}-\eqref{newA5} hold. If
	\begin{equation}\label{newA6}\tag{\textbf{A6}}
	\lim_{t\to+\infty}\sup_{x\in E}\ \underset{y\in E}{\mathrm{essup}}\ |\widetilde{p}(t,x,y)-1|=0,
	\end{equation}
	then, for any $\mu\in\me$ and $f\in \mathcal{B}^{+}(E)$ with $f/h$ bounded,
	$$\lim_{t\to+\infty}\e^{\lambda_{1}t}\langle f,X_{t}\rangle =(f,\hat{h}) W^{h}_{\infty}(X)\qquad\mbox{ in }L^{p}(\p_{\mu}).$$
\end{theorem}

\begin{theorem}[Strong law of large numbers]\label{them:slln}
	Suppose \eqref{newA1}-\eqref{newA6} hold.
	If
	\begin{equation}\label{newA7}\tag{\textbf{A7}}
	\lim_{t\to 0+}\|\widetilde{P}_{t}\phi-\phi\|_{\infty}=0\qquad\forall \phi\in C_{0}(E),
	\end{equation}
	where $C_{0}(E)$ denotes the space of bounded continuous functions that vanish at $\partial$, then, there exists $\Omega_{0}$ of $\p_{\mu}$-full probability for every $\mu\in\me$, such that on $\Omega_{0}$, for every $m$-almost everywhere continuous function $f$ with $f/h$ bounded,
	$$\lim_{t\to+\infty}\e^{\lambda_{1}t}\langle f,X_{t}\rangle =(f,\hat{h}) W^{h}_{\infty}(X).$$
\end{theorem}
The proofs of the above results will be given in Section \ref{sec:proofs}.

\begin{remark}
\rm
In this paper we assume that the spatial motion is a symmetric Borel right process. This assumption is not
necessary. An extension is possible, at least, to some extent.
One direction is to assume
that the spatial motion
is a transient Borel standard process on a Luzin space, which has a strong dual process. Definitions of smooth measures and the Kato class
can then be extended, while still preserving the properties used in this paper. We refer the readers to \cite{Chen,CS} for the Kato class
measures defined in this way.
Therefore, methods in this paper can be applied to establish the LLN for such superprocesses.
Nevertheless, we keep to the less general class of spatial motions in this paper for the sake of mathematical convenience.
\end{remark}

\subsection{Examples}\label{sec:examples}

In what follows, we will illustrate our main results by several concrete examples for which the branching mechanisms are local or non-local.

\begin{example}\label{eg:localbm}
	\rm
	In the case of a purely local branching mechanism where $\psi=\phi^{L}$ is given by \eqref{local part},
the auxiliary process $\widetilde{\xi}$ moves as a copy of the Doob $h$-transformed process $\xi^{h}$ of the spatial motion (cf. Proposition \ref{prop1} and Remark \ref{remark:many to one} below).
Therefore,
conditions \eqref{newA5}-\eqref{newA7}
are
reduced to the following:
	
	\noindent(\textbf{A5'})\hspace{0.5cm} $\sup_{x\in E}h(x)^{p-1}\int_{0,+\infty}u^{p}\pi^{L}(x,\ud u)<+\infty$ for some $p\in (1,2]$;
	
	\noindent(\textbf{A6'}) \hspace{.5cm} $\lim_{t\to +\infty}\sup_{x\in E}\mathrm{essup}_{y\in E}|p^{h}(t,x,y)-1|=0$;
	
	\noindent(\textbf{A7'}) \hspace{.5cm} $\lim_{t\to 0+}\|P^{h}_{t}\phi-\phi\|_{\infty}=0\quad \forall \phi\in C_{0}(E)$,\\
	where $P^{h}_{t}$ denotes the transition semigroup of
$\xi^{h}$
 and $p^{h}(t,x,y)$ denotes its transition density with respect to the measure $\widetilde{m}(\ud y)=h(y)^{2}m(\ud y)$.
	There is a large class of (purely local branching) superprocesses that satisfies conditions \eqref{newA1}-\eqref{newA7}, see for example, \cite[Examples 1,2,4,5]{CRY}. Therefore, Theorem \ref{them:wlln} and Theorem \ref{them:slln} can be applied to these superprocesses.
\end{example}

\begin{example}
	\label{eg:KP}
	\rm
		Suppose $E=\{1,2,\ldots,K\}$, $m$ is the counting measure on $E$ and $P_{t}f(i)= f(i)$ for all $i\in E$, $t\ge 0$ and $f\in\mathcal{B}^{+}(E)$. For $i\in E$ and $\mathbf{u}=(u_{1},u_{2},\cdots,u_{K})^{\mathrm{T}}\in [0,+\infty)^{K}$, define the function
		\begin{equation}
		\label{multitype}
		\psi(i,\boldsymbol{u}):=a_{i}u_{i}+b_{i}u_{i}^{2}- \boldsymbol{u}\cdot \boldsymbol{\eta}_i +\int_{(0,+\infty)^{K}}\left(\e^{-\boldsymbol{u}\cdot \boldsymbol{y}}-1+\boldsymbol{u}\cdot\boldsymbol{y}\right)\Gamma_{i}(\ud \boldsymbol{y}),
		\end{equation}
		where $\boldsymbol{u}\cdot \boldsymbol{y}=\sum_{i\in E} u_iy_i$ is the inner product of two vectors, $a_{i}\in (-\infty,+\infty)$, $b_{i}\ge 0$, $\boldsymbol{\eta}_i=(\eta_{i1},\cdots,\eta_{iK})^{\mathrm{T}}\in [0,\infty)^K$, and   $ \Gamma_{i}(\ud \boldsymbol{y})$ is a measure on $(0,+\infty)^{K}$ such that
		$$\int_{(0,+\infty)^{K}}(\boldsymbol{1}\cdot\boldsymbol{y})\wedge (\boldsymbol{1}\cdot\boldsymbol{y})^2 \Gamma_i(\ud  \boldsymbol{y})<+\infty\ \mbox{ and }\ \int_{(0,+\infty)^{K}}y_{j}\Gamma_{i}(\ud \boldsymbol{y})\le \eta_{ij}\quad\mbox{ for }i\not=j\in E.$$
		Without loss of generality
		we can assume that $\eta_{ii}=0$ for all $i\in E$
		(otherwise, we can
		change the value to $a_i$). We assume that there is a $p\in (1,2]$ such that
		\begin{equation}\label{eg1.condi}
		\int_{(0,+\infty)^{K}}\sum_{j=1}^{K}y_{j}^p\,\Gamma_{i}(\ud \boldsymbol{y})<+\infty,\qquad\forall i\in E.
		\end{equation}
		As a special case of the model given in Section \ref{sec:model}, we have a non-local branching superprocess $\{X_{t}:t\ge 0\}$ in $\me$ with transition probabilities given by
		\begin{equation}\nonumber
		\p_{\mu}\left[\exp\left(-\langle f,X_{t}\rangle\right)\right]=\exp\left(-\langle V_{t}f,\mu\rangle\right)\quad\mbox{ for }\mu\in\me,\ t\ge 0\mbox{ and }f\in \mathcal{B}^{+}_{b}(E),
		\end{equation}
		where $V_{t}f(i)$ is the unique nonnegative locally bounded solution to
		\begin{equation}\nonumber
		V_{t}f(i)=f(i)-\int_{0}^{t}\psi(i,V_{s}f)\ud s\qquad\mbox{ for } t\ge 0,\ i\in E.
		\end{equation}
		For every $i\in E$ and $\mu\in\me$, we define $\mu^{(i)}:=\mu(\{i\})$. The map \ $\mu\mapsto (\mu^{(1)},\cdots,\mu^{(K)})^{\mathrm{T}}$ \ is a homeomorphism between $\me$ and $[0,+\infty)^{K}$. Hence \ $\{(X^{(1)}_{t},\cdots,X^{(K)}_{t})^{\mathrm{T}}:t\ge 0\}$ \ is a Markov process in $[0,+\infty)^{K}$, which is called a
		{\it $K$-type continuous-state branching process}
($K$-type CSBP in abbreviation).
		Define the $K\times K$ matrix \ $\boldsymbol{M}(t)=(M(t)_{ij})_{i,j\in E}$ \ by \ $M(t)_{ij}:=\p_{\delta_{i}}\left[X^{(j)}_{t}\right]$ \ for $i,j\in E$.
		Let $\mathfrak{P}_{t}$
		denote the mean semigroup of $X$, that is
		$$
		\mathfrak{P}_{t}f(i)
		:=\p_{\delta_{i}}\left[\langle f,X_{t}\rangle\right]=\sum_{j=1}^{K}M(t)_{ij}f(j)\qquad\mbox{ for }i\in E,\ t\ge 0\mbox{ and }f\in\mathcal{B}^{+}(E).$$
		According to \cite[lemma 3.4]{blp1},
		\begin{equation}\nonumber
		\boldsymbol{M}(t)= {\rm e}^ {t\boldsymbol{A}^\mathrm{T}},\qquad t\geq 0,
		\end{equation}
		where the matrix $\boldsymbol{A}=(A_{ij})_{i,j}$  is given by $A_{ij}=-a_i\delta_{ij}+\eta_{ij}$ and $\boldsymbol{A}^\mathrm{T}$ is its transpose. If $\boldsymbol{A}^{\mathrm{T}}$ is irreducible, then Perron-Frobenius theory implies that there exist $\Lambda\in\R$ and right and left eigenvectors \ $\boldsymbol{h}, \hat{\boldsymbol{h}}\in\R_+^K$ \
with all their coordinates strictly positive
		such that
		$$\boldsymbol{M}(t)\boldsymbol{h} = {\rm e}^{\Lambda t}\boldsymbol{h} \qquad \mbox{ and }\qquad  \hat{\boldsymbol{h}}^{\mathrm{T}}\boldsymbol{M}(t)  = {\rm e}^{\Lambda t}\hat{\boldsymbol{h}}^{\mathrm{T}}, \qquad \mbox{ fot all }t\geq 0.$$
		For convenience we shall normalise $\boldsymbol{h}$ and $\hat{\boldsymbol{h}}$ such that \  $\boldsymbol{h}\cdot\boldsymbol{h}=\boldsymbol{h}\cdot\hat{\boldsymbol{h}}=1$. \
		Moreover, we have
		\begin{equation}
		e^{-\Lambda t}M(t)_{ij}\to h_{i}\hat{h}_{j}\quad\mbox{ as }t\to +\infty\quad\forall i,j\in E.\label{eg1}
		\end{equation}
		When $\Lambda\leq 0$, the $K$-type CSBP is extinct a.s.,
in other words, $$\p_{\mu}\left(\lim_{t\to+\infty}X^{(i)}_{t}=0\right)=1$$ for any $\mu\in\me$ (cf.
		\cite[Theorem 2]{KP} and \cite[Example 7.1]{RSY}).
		Henceforth we assume $\Lambda>0$.
		In view of \eqref{eg1.condi} and \eqref{eg1}, one can easily verify that conditions \eqref{newA1}-\eqref{newA5} hold for $\lambda_{1}=-\Lambda$, $\boldsymbol{h}$ and $\widehat{\boldsymbol{h}}$ and that \eqref{newA6} holds with \ $\widetilde{p}(t,i,j)=e^{-\Lambda t}(h_{i}\widehat{h}_{j})^{-1}M(t)_{ij}$ \ (cf.
		\cite[Example 7.1]{RSY}). The auxiliary process $\widetilde{\xi}$ is a finite-state Markov chain and hence is a Feller process. So condition \eqref{newA7} is automatically true.
We note that $W_{t}(X):=e^{-\Lambda t}\sum_{i=1}^{K}h_{i}X^{(i)}_{t}$ is a nonnegative martingale. Applying Theorems \ref{them:wlln} and \ref{them:slln}, we conclude that for every $\mu\in\me$ and $i\in E$, \ $\lim_{t\to+\infty}\e^{-\Lambda t}X^{(i)}_{t}=\hat{h}_{i}W_{\infty}$\  $\p_{\mu}$-a.s. and in $L^{p}(\p_{\mu})$, where $W_{\infty}$ denotes the martingale limit of $W_{t}(X)$. In particular on the event $\{W_{\infty}>0\}$,
		$$\lim_{t\to+\infty}\frac{X^{(i)}_{t}}{\sum_{j=1}^{K}X^{(j)}_{t}}=\frac{\hat{h}_{i}}{\sum_{j=1}^{K}\hat{h}_{j}}\quad \p_{\mu}\mbox{-a.s.}$$
		The a.s. convergence of this result is also obtained in \cite[Thoeorem 1.4]{KPR}.
	
\end{example}

\begin{example}\label{eg:multitype superprocess}\rm
	Suppose that $S=\{1,\ldots,K\}$, $D$ is a bounded $C^{1,1}$ domain in $\R^{d}$, and $m$ is the counting measure times the Lebesgue measure on $E=S\times D$.
	Suppose $\mathcal{L}^{(i)}$ ($i\in S$) is a second order differential operator of the form
	$$\mathcal{L}^{(i)}=\sum_{n,m=1}^{d}\frac{\partial}{\partial x_{n}}\left(\alpha^{(i)}_{n,m}(x)\frac{\partial}{\partial x_{m}}\right)\quad \mbox{ on }\R^{d},$$
	where for every $x\in\R^{d}$,  \ $\boldsymbol{B}^{(i)}(x):=\left(\alpha^{(i)}_{n,m}(x)\right)_{1\le n,m\leq d}$ \ is a uniformly elliptic symmetric matrix with \ $\alpha^{(i)}_{n,m}(x)\in C^{2,\gamma}(\R^{d})$ \ and $\gamma\in (0,1)$. Here, $C^{2,\gamma}(\R^{d})$ denotes the space of two times continuous differentiable functions whose second order derivatives are $\gamma$-H\"{o}lder continuous. Let $(\xi^{(i)},\Pi^{(i)})$ be a symmetric diffusion on $\R^{d}$ with generator $\mathcal{L}^{(i)}$, and $\xi^{(i),D}$ be the subprocess of $\xi^{(i)}$ killed upon leaving $D$. It is known that the semigroup of $\xi^{(i),D}$, denoted by $P^{(i),D}_{t}$, admits a transition density function \ $p^{(i)}_{D}(t,x,y)$ \ with respect to the Lebesgue measure, which is symmetric in $(x,y)$ for each $t>0$. For $f\in \mathcal{B}^{+}(E)$, we use the convention $\boldsymbol{f}(x)=(f(1,x),\ldots,f(K,x))^{\mathrm{T}}=(f_{1}(x),\ldots,f_{K}(x))^{\mathrm{T}}$. Let $\xi$ be a Markov process on $E$ with semigroup \  $P_{t}f(i,x):=P^{(i),D}_{t}f_{i}(x)$ for $f\in\mathcal{B}^{+}(E)$ and $(i,x)\in E$.
	Define the branching mechanism
	$$\Psi((i,x),\boldsymbol{f}):=\psi(i, \boldsymbol{f}(x)),\qquad \forall (i,x)\in E,\ f\in\mathcal{B}^{+}(E),$$
	where $\psi(i,\cdot)$ is given by
\eqref{multitype}-\eqref{eg1.condi}
in Example \ref{eg:KP}.
	Suppose \ $\{X_{t}:t\ge 0\}$ \ is a $(P_{t},\Psi)$-superprocess in $\me$.
	For every $i\in E$ and $\mu\in\me$, we define \ $\mu^{(i)}(A):=\mu(\{i\}\times A)$.\  The map \ $\mu\mapsto (\mu^{(1)},\cdots,\mu^{(K)})^{\mathrm{T}}$ \ is a homeomorphism between $\me$ and $\mathcal{M}(D)^{K}$. Hence, \ $\{(X^{(1)}_{t},\cdots,X^{(K)}_{t})^{\mathrm{T}}:t\ge 0\}$ is a Markov process in $\mathcal{M}(D)^{K}$, which is called a {\it $K$-type superdiffusion}.
	
	Let denote by $\mathfrak{P}_{t}$  the mean semigroup of $X_{t}$, that is,
	$$\mathfrak{P}_{t}f(i,x)=\mathrm{P}_{\delta_{(i,x)}}\left[\sum_{j=1}^{K}\langle f_{j},X^{(j)}_{t}\rangle\right].$$
In view of \eqref{lem5.1.1} below, we have
	\begin{equation}\nonumber
	\mathfrak{P}_{t}f(i,x)=\e^{-a_{i}t}P^{(i),D}_{t}f_{i}(x)+\int_{0}^{t}\e^{-a_{i}s}\sum_{j\not= i}\eta_{ij}P^{(j),D}_{s}(\mathfrak{P}_{t-s}(j,\cdot))(x)\ud s
	\end{equation}
	for every $f\in \mathcal{B}_{b}(E)$ and $(i,x)\in E$, where  $a_i$ and $\eta_{ij}$ are the linear local and non-local parts
	of $\psi(i,\cdot)$, respectively. Now let $a_{0}:=\max_{i\in S}\left(-a_{i}+\sum_{j\in S}\eta_{ij}\right)+1$ and $Q_{t}f(i,x):=e^{-a_{0} t}\mathfrak{P}_{t}f(i,x)$ \ for all $f\in\mathcal{B}_{b}(E)$. Then, $Q_{t}$ satisfies that
	\begin{eqnarray}\nonumber
	Q_{t}f(i,x)&=&\e^{-(a_{i}+a_{0})t}P^{(i),D}_{t}f_{i}(x)\nonumber\\
	&&+\int_{0}^{t}\e^{-(a_{i}+a_{0})s}\left(a_{i}+a_{0}\right)\sum_{j\not= i}\frac{\eta_{ij}}{a_i+a_0}P^{(j),D}_{s}(Q_{t-s}(j,\cdot))(x)\ud s.\nonumber
	\end{eqnarray}
	This implies that $Q_{t}$ is the semigroup of
	a switched diffusion $(\Theta_{t},\Xi_{t})_{t\ge 0}$
	on $E$ with generator
	\[
	\mathcal{S}u=\begin{pmatrix}
	\mathcal{L}^{(1)}&0&\cdots&0\\
	0&\mathcal{L}^{(2)}&\cdots&0\\
	\vdots&\vdots&\ddots&\vdots\\
	0&0&\cdots&\mathcal{L}^{(K)}
	\end{pmatrix}u+\boldsymbol{Q}u
	, \qquad\qquad\forall u:D\mapsto \R^{K}
	\]
	where \ $\boldsymbol{Q}=(q_{ij})_{1\le i,j\le K}$ \ is given by \ $-q_{ii}=a_{i}+a_{0}$ \ and \ $q_{ij}=\eta_{ij}$ \ for $i\not= j$.
	The movement of the switched diffusion $(\Theta_{t},\Xi_{t})_{t\ge 0}$ is described as follows: The process $\Theta$ moves as
an $S$-valued Markov chain with intensity matrix $\boldsymbol{Q}$. When $\Theta$ is in a state $j\in S$, $\Xi$ moves as an independent copy of $\xi^{(j),D}$ until
$\Theta$ has a jump.
When $\Theta$ changes from $j$ to another state $k\in S$, $\Xi$
immediately and continuously evolves
as an independent copy of $\xi^{(k),D}$ and so on, until the lifetime of $\Theta$.
	It is known by \cite[Theorem 5.3]{CZ} that $Q_{t}$ admits a transition density function \ $q(t,(i,x),(j,y))$ \ with respect to $m$. Moreover by \cite[Theorem 5.3]{CZ} one can deduce that for every $i,j\in S$ and $t>0$, \ $(x,y)\mapsto q(t,(i,x),(j,y))$ is jointly continuous and
	that there are positive constants $t_{0},\ c_{i},i=1,\ldots,4$ such that
	\begin{equation}\label{eg6.3.1}
	c_{1}p_{0}(c_{2}t,x,y)\le q(t,(i,x),(j,y))\le c_{3}p_{0}(c_{4}t,x,y),\qquad \forall t\in (0,t_{0}], \ i,j\in S,
	\end{equation}
	where  \ $p_{0}(t,x,y)$  is the transition density of a killed Brownian motion in $D$. It follows imme\-dia\-tely that $\mathfrak{P}_{t}$ has a transition density \ $\mathfrak{p}(t,(i,x),(j,y))$ with respect to $m$ which satisfies that
	\begin{equation}\nonumber
	c_{1}\e^{a_{0}}p_{0}(c_{2}t,x,y)\le \mathfrak{p}(t,(i,x),(j,y))\le c_{3}\e^{a_{0} t}p_{0}(c_{4}t,x,y),\qquad \forall t\in (0,t_{0}],\  i,j\in S.
	\end{equation}
	Note that \ $\int_{E^{2}}\mathfrak{p}(t,(i,x),(j,y))^{2}m(\ud i,\ud x)m(\ud j,\ud y)<+\infty$ \ for every $t\in (0,t_{0}]$. Thus $\mathfrak{P}_{t}$ is a Hilbert-Schmidt operator in $L^{2}(E,m)$ and hence is compact. The same is true for its dual operator $\hat{\mathfrak{P}}_{t}$. If we use $\sigma(\mathcal{L)}$ and $\sigma(\hat{\mathcal{L}})$ to denote the spectrum of the generators of $\mathfrak{P}_{t}$ and $\hat{\mathfrak{P}}_{t}$ respectively, then it follows by Jentzch's theorem that \ $-\lambda_{1}:=\sup \Re (\sigma(\mathcal{L}))=\sup \Re (\sigma(\hat{\mathcal{L}}))$ \ is a simple eigenvalue for both $\mathcal{L}$ and $\hat{\mathcal{L}}$, and that an eigenfunction $h$ of $\mathcal{L}$ associated with $-\lambda_{1}$ and an eigenfunction $\hat{h}$ of $\hat{\mathcal{L}}$ associated with $-\lambda_{1}$ can be chosen strictly positive on $E$ and
	satisfying
	$\int_{E}h^{2}\ud m=\int_{E}h\hat{h}\ud m=1$. \
	It is proved in \cite[Section 3]{CRS} that there exists a constant $c_{5}>1$ such that
	\begin{equation}\label{eg6.3.2}
	c^{-1}_{5}\delta_{D}(x)\le h(i,x),\hat{h}(i,x)\le c_{5}\delta_{D}(x),\qquad \forall (i,x)\in E.
	\end{equation}
	Here, $\delta_{D}(x)$ denotes the Euclidean distance between $x$ and the boundary of $D$. We assume $\lambda_{1}<0$. One can easily verify that conditions \eqref{newA1}-\eqref{newA5} hold for this example.
Moreover,
	\begin{equation}\label{eg:6.3.3} \widetilde{p}(t,(i,x),(j,y))=\frac{\e^{\lambda_{1}t}\mathfrak{p}(t,(i,x),(j,y))}{h(i,x)\hat{h}(j,y)}=\frac{\e^{(\lambda_{1}+a_{0})t}q(t,(i,x),(j,y))}{h(i,x)\hat{h}(j,y)}.
	\end{equation}
	In view of \eqref{eg6.3.1} and \eqref{eg6.3.2}, one can apply
a similar argument
as in \cite[section 3]{CRS} to show that the semigroup $Q_{t}$ is intrinsically ultracontractive, i.e. for any $t>0$ there is a constant $c_{t}>0$ such that
	$$q(t,(i,x),(j,y))\le c_{t}h(i,x)\hat{h}(j,y),\qquad\forall (i,x),(i,y)\in E.$$
	As a consequence, there exist constants $t_{1},c_{6},c_{7}>0$ such that
	$$\left|\widetilde{p}(t,(i,x),(j,y))-1\right|\le c_{6}\e^{-c_{7}t},\qquad \forall t>t_{1},\ (i,x),(j,y)\in E.$$
	Hence, condition \eqref{newA6} is satisfied. It is known that there are positive constants $C_{i},i=1,\cdots,4$ such that for any $t\in (0,1]$ and $x,y\in D$,
	\begin{multline}\label{eg:6.3.4}
	C_1\left(\frac{\delta_{D}(x)}{\sqrt{t}}\wedge 1\right)\left(\frac{\delta_{D}(y)}{\sqrt{t}}\wedge 1\right)t^{-\frac{d}{2}}\e^{-\frac{C_2|x-y|^{2}}{2}}\le p_{0}(t,x,y)\\
	\le C_3\left(\frac{\delta_{D}(x)}{\sqrt{t}}\wedge 1\right)\left(\frac{\delta_{D}(y)}{\sqrt{t}}\wedge 1\right)t^{-\frac{d}{2}}\e^{-\frac{C_4|x-y|^{2}}{2}}.
	\end{multline}
	In view of \eqref{eg6.3.1}-\eqref{eg:6.3.4}, we can apply similar calculations as in \cite[Example 3]{CRY} to show condition \eqref{newA7} is satisfied. Let  $W_{\infty}$ denote the limit of the nonnegative martingale \ $W_{t}=e^{\lambda_{1} t}\sum_{j=1}^{K}\langle h_{j},X^{(j)}_{t}\rangle$. \
	Applying Theorems \ref{them:wlln} and \ref{them:slln}, we conclude that for every $\mu\in\me$, $i\in S$ and $f_{i}\in C^{+}(D)$ with $f_{i}/h_{i}$ bounded, \ $\lim_{t\to+\infty}\e^{\lambda_{1}t}\langle f_{i},X^{(i)}_{t}\rangle=W_{\infty}\int_{D}f_{i}(x)\hat{h}_{i}(x)\ud x$ \ $\p_{\mu}$-a.s. and in $L^{p}(\p_{\mu})$.
\end{example}

\begin{example}\label{eg:CKS}{\rm
		Suppose that $E$ is a bounded $C^{1, 1}$ open set in
		$\R^{d}\ (\ud \ge 1)$, $m$ is the Lebesgue
		measure on $E$, $\alpha\in (0, 2)$, $\beta\in [0, \alpha\wedge d)$ and that $\xi=(\xi_t, \Pi_x)$ is an $m$-symmetric Hunt process on $E$ satisfying the following conditions: (1) $\xi$ has a L\'evy system
		$(N,t)$
		where $N=N(x, \ud y)$ is a jumping kernel given by
		$$
		N(x, \ud y)=\frac{C_1}{|x-y|^{d+\alpha}}\,\ud y, \qquad x, y\in E
		$$
		for some constant $C_1>0$.
		
		\noindent(2)  $\xi$ admits a
		jointly continuous transition density $p(t, x, y)$ with respect to $m$
		and that there exists a constant $C_2>1$ such that
		$$
		C_2^{-1}q_{\beta}(t, x, y)\le p(t, x, y)\le C_2q_{\beta}(t, x, y),
		\qquad
		\forall
		(t, x, y)\in (0, 1]\times E\times E,
		$$
		where
		\begin{equation}\label{e:CKSassump}
		q_{\beta}(t, x, y)=\left(1\wedge\frac{\delta_E(x)}{t^{1/\alpha}}\right)^\beta
		\left(1\wedge\frac{\delta_E(y)}{t^{1/\alpha}}\right)^\beta
		\left(t^{-d/\alpha}\wedge\frac{t}{|x-y|^{d+\alpha}}\right).
		\end{equation}
		Here $\delta_E(x)$ stands for
		the Euclidean distance between $x$ and the boundary of $E$.
One concrete example of $\xi$ is the killed symmetric $\alpha$-stable process
in $E$. In this case, \eqref{e:CKSassump} is satisfied with $\beta=\alpha/2$.
Another example of $\xi$
is the censored symmetric $\alpha$-stable process
in $E$ introduced in \cite{BBC} when $\alpha\in (1, 2)$. In this case, \eqref{e:CKSassump} is satisfied with $\beta=\alpha-1$.
In fact, by using \cite{CKS}, one could also include the case when $E$ is a $d$-set,
$\alpha\in (0, 2)$ and $\xi$ is an $\alpha$-stable-like process in $E$.

Suppose that the branching mechanism $\psi$ is given by \eqref{bm1}
		and satisfies \eqref{newA4}.
		We further assume that the kernel
		$\gamma(x, \ud y)$ has a density $\gamma(x, y)$ with respect to $m$, which satisfies that
		$$
		\gamma(x, y)\le C_3|x-y|^{\epsilon-d} \qquad \forall x, y\in E
		$$
		for some $C_3,\ \epsilon>0$. It is proved in \cite[Example 7.3]{RSY} that \eqref{newA1}-\eqref{newA3} and \eqref{newA6} are satisfied.
		They also proved that the mean semigroup $\mathfrak{P}_{t}$ of this superprocess has a density function $\mathfrak{p}(t,x,y)$ with respect to $m$ such that $(x,y)\mapsto \mathfrak{p}(t,x,y)$ is jointly continuous for each $t>0$ and
		\begin{equation}
		C_{4}^{-1}q_{\beta}(t,x,y)\le\mathfrak{p}(t,x,y)\le C_{4}q_{\beta}(t,x,y)\quad \forall t\in (0,1],\ x,y\in E,
		\end{equation}
		for some  $C_{4}>1$. Moreover,
		\begin{equation}\label{eg6.5.6}
		C^{-1}_{5}
		\delta_E(x)^{\beta}\le h(x),\hat{h}(x)\le
		C_{5}
		\delta_E(x)^{\beta}\qquad\forall x\in E
		\end{equation}
		for some $C_{5}>1$.
		In view of \eqref{prop2.0 density} and \eqref{e:CKSassump}-\eqref{eg6.5.6}, we can show condition \eqref{newA7} by applying similar calculations as in \cite[Example 5]{CRY}.
		Therefore, Theorems \ref{them:wlln} and \ref{them:slln} can be applied to this example as long as condition \eqref{newA5} holds. By \eqref{eg6.5.6}, condition \eqref{newA5} is satisfied if and only if  there is a $p\in (1,2]$ such that
		\begin{equation}
		\sup_{x\in E}\delta_{E}^{-\beta}(x)\int_{\me^{0}}\nu(\delta^{\beta}_{E})^{p}H(x,\ud \nu)<+\infty.\label{eg:6.5.1}
		\end{equation}

		An example of a branching mechanism that satisfies \eqref{eg:6.5.1} and cannot be
		decomposed
		into local and non-local parts is
		\begin{equation}\label{eq:non-decomp bm}
		\begin{split}
		\psi(x,f)=&a(x)f(x)+b(x)f(x)^{2}\\
		&+c(x)\int_{0}^{1}\left(\exp\{-u f(x)-u^{2}\pi(x,f)\}-1+u f(x)\right)\frac{1}{u^{1+\theta}}\ud u,
		\end{split}
		\end{equation}
		where \ $\theta\in (1,2)$, \ $a(x)\in C_{b}(E)$, \ $b(x),c(x)\in C^{+}_{b}(E)$, \ $c(x)\le C_{6}\delta_{E}(x)^{\beta}$ \  for some $C_{6}>0$, and $\pi(x,\ud y)$ is a probability kernel on $E$ which has a density function $\pi(x,y)$ with respect to $m$ satisfying that \ $\pi(x,y)\le C_{7}|x-y|^{\epsilon-d}$ \  for some $\epsilon,C_{7}>0$. In fact, $\psi$ given by \eqref{eq:non-decomp bm} can be represented in the form of \eqref{1.1} with $H(x,\ud \nu)$ being the image of \ $c(x)u^{-1-\theta}\ud u$ \ under the mapping \ $u\mapsto u\delta_{x}(\ud y)+u^{2}\pi(x,\ud y)$ \  of $[0,1]$ into $\me^{0}$, and \ $\gamma(x,\ud y)=\frac{c(x)}{2-\theta}\pi(x,\ud y)$.}
\end{example}

\section{Martingale problems and representation of superprocesses}\label{sec:martingale and representation}

The martingale problem of superprocesses with branching mechanisms given by \eqref{bm1} is studied in \cite{Li} under some Feller type assumptions and the assumption
$$\sup_{x\in E}\int_{\me^{0}}\nu(1)^{2}H(x,\ud \nu)<+\infty.$$
These conditions
guarantee that the martingale measure induced by the martingale problem is \textit{worthy}. Using the worthy martingale measure, \cite{Li} establishes a representation for superprocesses in terms of stochastic integrals.
In this section,
we shall drop the above
$L^{2}$-moment
condition and investigate the martingale problem under much weaker hypotheses. As a result, we obtain the same type of representation for superprocesses when the underlying martingale measures are not necessarily orthogonal or worthy.
All martingales or local martingales mentioned in this section will be relative to the filtration $(\mathcal{F}_{t})_{t\ge 0}$ and the probability $\p_{\mu}$ where $\mu\in\me$.

\subsection{Martingale problems of superprocesses}\label{sec:matingale problem}

A measurable function $f$ is said to be \textit{finely continuous} relative to  $\xi$ if \ $t\mapsto f(\xi_{t})$ \ is a.s. right continuous on $[0,+\infty)$.
Let $U^{\alpha}$ denote the $\alpha$-resolvent of $(P_{t})_{t\ge 0}$, (in other words \ $U^{\alpha}f(x):=\int_0^\infty e^{-\alpha t}P_{t}f(x)\ud t$). Recall that $\mathcal{C}^{\xi}_{b}(E)$ is the set of bounded measurable functions that are finely continuous with respect to $\xi$. Fix an arbitrary $\beta>0$, define \ $D(\mathcal{A}):=U^{\beta}\mathcal{C}^{\xi}_{b}(E)$, \  and for any \ $f=U^{\beta}g\in D(\mathcal{A})$ \ with \ $g\in \mathcal{C}^{\xi}_{b}(E)$,\  set \ $\A f:=\beta f-g$. It is known
(cf.
\cite[A.6]{Li}) that \ $(\A,D(\A)):D(\A)\to \C$ \ defines a linear operator which is independent of $\beta$. Moreover for every \ $f\in D(\A)$, \ $(P_{t}f-f)/t$ \ converges boundedly and pointwise to $\A f$ as $t\to 0$. We call \ $(\A,D(\A))$ \ the \textit{weak generator} of $\xi$.
For a measurable function $f$, we set
 $$e_{f}(t):=\exp\left(-\int_{0}^{t}f(\xi_{s})\ud s\right),\qquad\forall
 t\ge 0,
 $$
whenever it is well defined.

\begin{lemma}\label{lem3.1}
Suppose $f\in  \mathcal{B}_{b}(E)$. If $f$ is finely continuous with respect to $\xi$, then \ $t\mapsto \langle f,X_{t}\rangle$ is right continuous on $[0,+\infty)$ almost surely. If \ $t\mapsto f(\xi_{t})$ \ has left limits on $(0,+\infty)$ a.s., then so does \ $t\mapsto \langle f,X_{t}\rangle$.
\end{lemma}

\proof The idea of this proof is from \cite[Theorem 3.5(a)]{F1988}.
For a function $g\in\mathcal{B}_{b}(E)$, we use $(P^{g}_{t})_{t\ge 0}$ to denote the Feynman-Kac semigroup given by
\begin{equation}
P^{g}_{t}f(x)=\Pi_{x}\left[e_{g}(t)f(\xi_{t})\right],\qquad x\in E,\ f\in\mathcal{B}_{b}(E).\nonumber
\end{equation}
Using this notation, one can rewrite $\mathfrak{P}_{t}f(x)$ given in \eqref{mean equation} as follows
\begin{equation}
\mathfrak{P}_{t}f(x)=P^{a }_{t}f(x)+\int_{0}^{t}P^{a }_{t-s}\gamma (\mathfrak{P}_{s}f)(x)\ud s,\qquad t\ge 0,\ x\in E,\ f\in\mathcal{B}_{b}(E).\label{lem5.1.1}
\end{equation}
 By Gronwall's inequality, we have
\begin{equation}\label{mean ineq}
\|\mathfrak{P}_{t}f\|_{\infty}\le \e^{c_{0}t}\|f\|_{\infty},\qquad \forall t\ge 0,\ f\in \mathcal{B}_{b}(E),
\end{equation}
where \ $c_{0}:=\|\gamma(\cdot,1)\|_{\infty}+\|a^{-}\|_{\infty}$.
Now choose an arbitrary constant $q_{0}>c_{0}$. For $t\ge 0$, define the operator \ $Q_{t}:\mathcal{B}_{b}(E)\to \mathcal{B}_{b}(E)$ \ by \ $Q_{t}f(x):=\e^{-q_{0} t}\mathfrak{P}_{t}f(x)$. \ It follows by \eqref{lem5.1.1} that
\begin{eqnarray}
Q_{t}f(x)&=&P^{a +q_{0}}_{t}f(x)+\int_{0}^{t}P^{a +q_{0}}_{t-s}\gamma (Q_{s}f)(x)\ud s\nonumber\\
&=&P^{a +q_{0}}_{t}f(x)+\Pi_{x}\left[\int_{0}^{t}e_{a +q_{0}}(s)\left(a (\xi_{s})+q_{0}\right)\hat{\kappa}(\xi_{s},Q_{t-s}f)\ud s\right],\nonumber
\end{eqnarray}
where \ $\hat{\kappa}(x,\ud y):=\gamma (x,\ud y)/(a (x)+q_{0})$ \ is a sub-Markov kernel on $E$. We extend $\hat{\kappa}(x,\ud y)$ to a Markov kernel from $E$ to $E\cup\{\partial\}$ by setting \ $\hat{\kappa}(x,\{\partial\})=1-\int_{E}\hat{\kappa}(x,\ud y)$.
Let $\hat{\xi}$ be the Markov process obtained through a ``piecing out" procedure of \cite{INW}
(see also Section \ref{sec:interpretation})
from an infinite sequence of copies of the $e_{a+q_{0}}(t)$-subprocess of $\xi$, and the instantaneous distribution $\hat{\kappa}$.
Then $(Q_{t})_{t\ge 0}$ \ defined above is the semigroup corresponding to $\hat{\xi}$. It follows by \cite[Theorem A.43]{Li} that $(Q_{t})_{t\ge 0}$ \  induces the same fine topology on $E$ as $(P_{t})_{t\ge 0}$.

Now fix an arbitrary $\mu\in\me$. Let $\{T_{n}:n\ge 1\}$ be a decreasing sequence of bounded $\mathcal{F}_{t}$-stopping times with limit $T$. Define $\nu_{n}\in \me$ by
$$\nu_{n}(f):=\p_{\mu}\left[\e^{-q_{0} T_{n}}\langle f,X_{T_{n}}\rangle\right],\qquad \forall f\in \mathcal{B}_{b}(E),$$
and define $\nu$ analogously with \ $T_{n}$ \ replaced by $T$. Let \ $\mathfrak{U}^{q_{0}}f(x):=\int_{0}^{+\infty}Q_{t}f(x)\ud t$ \ for all $f\in\mathcal{B}_{b}(E)$.
One can easily show by strong Markov property and Fubini's theorem that
\begin{equation*}
\nu_{n}\left(\mathfrak{U}^{q_{0}}f\right)=\p_{\mu}\left[\int_{T_{n}}^{+\infty}\e^{-q_{0} s}\langle f,X_{s}\rangle \ud s\right]
\end{equation*}
for every $f\in\mathcal{B}_{b}(E)$. Hence, \ $\nu_{n}(\mathfrak{U}^{q_{0}}f)\ \uparrow \ \nu(\mathfrak{U}^{q_{0}}f)$ \ as $n\to +\infty$.
If $f$ is finely continuous relative to $\xi$, it is also finely continuous relative to $\hat{\xi}$, and so by \cite[Proposition 3.3]{F1988}, we have $\nu_{n}(f)\to \nu(f)$. \ Since $\{T_{n}:n\ge 1\}$ is arbitrary, \cite[VI.48]{DM1} yields the almost sure right continuity of \ $t\mapsto \e^{-q_{0} t}\langle f,X_{t}\rangle$. \
Hence,
we prove the first assertion. The second assertion follows analogously from \cite[Proposition 3.4(a)]{F1988}.\qed

\bigskip

We note that by definition every $f\in D(\A)$ is a $\beta$-excessive function relative to $(P_{t})_{t\ge 0}$ and thus \ $t\mapsto f(\xi_{t})$ \ is c\`{a}dl\`{a}g almost surely. This together with Lemma \ref{lem3.1} implies that \ $t\mapsto \langle f,X_{t}\rangle$ \ is c\`{a}dl\`{a}g almost surely for every $f\in D(\A)$.

 Let $N(\ud s,\ud \nu)$ be the random measure on $\R^{+}\times \mathcal{M}(E)^{0}$ defined by
\begin{equation*}
N(\ud s,\ud \nu):=\sum_{s\ge 0}\1_{\{\Delta X_{s}\not=0\}}\delta_{(s,\Delta X_{s})}(\ud s,\ud \nu).
\end{equation*}
Here, we use the standard notation \ $\Delta X_{s}:=X_{s}-X_{s-}$ \ for the jump of $X$ at time $s$. Let $\hat{N}(\ud s,\ud \nu)$ be the predictable compensator of $N(\ud s,\ud \nu)$ and $\widetilde{N}(\ud s,\ud \nu):=N(\ud s,\ud \nu)-\hat{N}(\ud s,\ud \nu)$ be the compensated random measure.
In view of condition \eqref{newA4} and the argument above, one can prove the following result in the same way as \cite[Theorem 7.13]{Li}.
\begin{theorem}\label{them:martingale}
Suppose \eqref{newA4} holds. The following statements are true.
\begin{enumerate}
\item [(i)] The predictable compensator $\hat{N}(\ud s,\ud \nu)$ is given by
$$\hat{N}(\ud s,\ud \nu)=\ud s\int_{E}X_{s-}(\ud x)H(x,\ud \nu)$$
where $H(x,\ud \nu)$ is the kernel associated with the non-linear part in \eqref{bm1}.
\item[(ii)] The superprocess $(X_{t})_{t\ge 0}$ has no negative jumps,
that is, for any $s\ge 0$, $\Delta X_{s}$ is a nonnegative random measure.
For any $f\in D(\A)$, the process
$$M_t(f):=\langle f,X_{t}\rangle -\langle f,X_{0}\rangle-\int_{0}^{t}\langle \A f+\gamma f-a f,X_{s}\rangle \ud s,\qquad t\ge 0,$$
is a
c\`{a}dl\`{a}g
martingale.
\item[(iii)]  For all $f\in D(\A)$,
$M_{t}(f)$
has a unique decomposition
$$M_t(f) =M^{c}_{t}(f)+M^{d}_{t}(f),\qquad t\ge 0,$$
where \ $t\mapsto M^{c}_{t}(f)$ \ is a square integrable continuous martingale with quadratic variation \ $\langle M^{c}(f)\rangle_{t}=2\int_{0}^{t}\langle bf^{2},X_{s}\rangle \ud s$, \ and
$$t \mapsto M^{d}_{t}(f)=\int_{0}^{t}\int_{\mathcal{M}(E)^{0}}\nu(f)\widetilde{N}(\ud s,\ud \nu)$$
is a purely discontinuous martingale.
\end{enumerate}
\end{theorem}

\subsection{A representation for superprocesses}\label{sec:representation}

For $f\in D(\A)$, it will be convenient to write
$$M_{t}(f)=\int_{0}^{t}\int_{E}\1_{\{s\le t\}}f(x)M(\ud s,\ud x).$$
We shall show in the following that the stochastic integral \ $\int_{0}^{t}\int_{E}\varphi(s,x)M(\ud s,\ud x)$ \ can be defined formally for a
large
class of integrands
$\varphi(s,x)$, which includes the functions $\{1_{\{s\le t\}}f(x):t\ge 0,f\in D(\A)\}$ as a subclass.

Let $(\omega,s,\nu)\mapsto F(\omega,s,\nu)$ be a predictable function on $\mathcal{W}^{+}_{0}\times \R^{+}\times \me^{0}$ such that
\begin{equation}\label{condi:integral N}
\p_{\mu}\left[\left(\sum_{s\le t}F(s,\Delta X_{s})^{2}\1_{\{\Delta X_{s}\not=0\}}\right)^{1/2}\right]<+\infty,\qquad\forall t\ge 0.
\end{equation}
Then, following \cite[Section II.1d]{Jacod}, one can define the stochastic integral of $F$ with respect to the compensated measure $\widetilde{N}(\ud s,\ud \nu)$, denoted by
$$\int_{0}^{t}\int_{\me^{0}}F(s,\nu)\widetilde{N}(\ud s,\ud \nu),$$
as the unique purely discontinuous local martingale whose jumps are indistinguishable from the process \ $F(s,\Delta X_{s})\1_{\{\Delta X_{s}\not=0\}}$. \ Condition \eqref{condi:integral N} holds in the special case where $F(\omega,s,\nu)=F_{\varphi}(\omega,s,\nu)=\int_{E}\varphi(s,x)\nu(\ud x)$
and $\varphi$ is a bounded measurable function on $\R^{+}\times E$.
Indeed, in this case we have
\begin{eqnarray}
\p_{\mu}\hspace{-.1cm}\left[\left(\sum_{s\le t}F_{\varphi}(s,\Delta X_{s})^{2}\1_{\{\Delta X_{s}\not=0\}}\hspace{-.1cm}\right)^{1/2}\right]\nonumber
\le \underset{x\in E}{\sup_{s\geq0}}\varphi(s,x)|\p_{\mu}\hspace{-.1cm}\left[\left(\sum_{s\le t}\Delta X_{s}(1)^{2}\1_{\{\Delta X_{s}\not=0\}}\hspace{-.1cm}\right)^{1/2}\right]\nonumber
\end{eqnarray}
where
$\Delta X_{s}(1)=\langle 1,\Delta X_{s}\rangle$. Moreover, we have
\begin{align}
 &\p_{\mu}\left[\left(\sum_{s\le t}\Delta X_{s}(1)^{2}\1_{\{\Delta X_{s}\not=0\}}\right)^{1/2}\right]\nonumber\\
&\le\p_{\mu}\left[\left(\sum_{s\le t}\Delta X_{s}(1)^{2}\1_{\{\Delta X_{s}(1)\le 1\}}\right)^{1/2}\right]
+\p_{\mu}\left[\left(\sum_{s\le t}\Delta X_{s}(1)^{2}\1_{\{\Delta X_{s}(1)>1\}}\right)^{1/2}\right]\nonumber\\
&\le\p_{\mu}\left[\sum_{s\le t}\Delta X_{s}(1)^{2}\1_{\{\Delta X_{s}(1)\le 1\}}\right]^{1/2}
+\p_{\mu}\left[\sum_{s\le t}\Delta X_{s}(1)\1_{\{\Delta X_{s}(1)>1\}}\right]\nonumber
\end{align}
In the first inequality, we use the fact that $(a+b)^{1/2}\le a^{1/2}+b^{1/2}$ for any $a,b\ge 0$. In the second inequality, we use Jensen's inequality and the fact that $(a_{1}+\cdots+a_{n})^{1/2}\le a_{1}^{1/2}+\cdots+a_{n}^{1/2}$ for any $n\ge 1$ and $a_{1},\ldots,a_{n}\ge 0$, respectively, to get the first and second term.
Therefore, by Theorem \ref{them:martingale}(i) we get
\begin{multline}
\p_{\mu}\left[\left(\sum_{s\le t}\Delta X_{s}(1)^{2}\1_{\{\Delta X_{s}\not=0\}}\right)^{1/2}\right]\\
\leq\p_{\mu}\left[\int_{0}^{t}\ud s\int_{E}X_{s-}(\ud x)\int_{\me^{0}}\nu(1)^{2}\1_{\{\nu(1)\le 1\}}H(x,\ud \nu)\right]^{1/2}\\
\quad+\p_{\mu}\left[\int_{0}^{t}\ud s\int_{E}X_{s-}(\ud x)\int_{\me^{0}}\nu(1)\1_{\{\nu(1)>1\}}H(x,\ud \nu)\right].\label{eq3.15}
\end{multline}
In view of the fact that $\nu(1)\wedge \nu(1)^{2}H(x,\ud \nu)$ is a bounded kernel from $E$ to $\me^{0}$, we can show by \eqref{mean ineq} that the expectations
on the right hand side
are finite.
In the sequel, we will write
$$M^{d}_{t}(\varphi)=\int_{0}^{t}\int_{E}\varphi(s,x)M^{d}(\ud s,\ud x):=\int_{0}^{t}\int_{\me^{0}}F_{\varphi}(s,\nu)\widetilde{N}(\ud s,\ud \nu)$$
for every $\varphi\in \mathcal{B}_{b}(\R^{+}\times E)$.

Define a random measure $\eta$ on $\R^{+}\times E\times E$ by
$$\eta(\ud s,\ud x,\ud y):= \ud s \int_{E}X_{s}(\ud z)2 b(z)\delta_{z}(\ud x)\delta_{z}(\ud y).$$
Immediately by \eqref{mean ineq} we have
\begin{equation*}
\p_{\mu}\left[\left|\int_{0}^{t}\int_{E^{2}}\eta(\ud s,\ud x,\ud y)\right|\right]\le 2\|b\|_{\infty}\p_{\mu}\left[\int_{0}^{t}X_{s}(1)\ud s\right]<+\infty.
\end{equation*}
Theorem \ref{them:martingale}(iii) yields that \ $\langle M^{c}(f)\rangle_{t} =\int_{0}^{t}\int_{E^{2}}f(x)f(y)\eta(\ud s,\ud x,\ud y)$ \ for every $f\in D(\A)$. Thus, by Doob's martingale inequality
\begin{equation*}
\p_{\mu}\left[\sup_{0\le s\le t}|M^{c}_{s}(f)-M^{c}_{s}(g)|^{2}\right]\le 4\p_{\mu}\left[\int_{0}^{t}\int_{E^{2}}|f(x)-g(x)||f(y)-g(y)|\eta(\ud s,\ud x,\ud y)\right]
\end{equation*}
for all $f,g\in D(\A)$. Using the above two inequalities and the fact that any element of $C_{b}(E)$ is the bounded pointwise limit of a sequence from $D(\A)$, one can extend the linear map $D(\A)\ni f\mapsto M^{c}(f)$ to a martingale functional \ $\{M^{c}(f),f\in C_{b}(E)\}$ \
in the same way as
\cite[Section 7.3]{Li}, and then further extend it to a martingale measure \ $M^{c}(\ud s,\ud x)$ \ on $\R^{+}\times E$, which satisfies that
$$M^{c}_{t}(f)=\int_{0}^{t}\int_{E}f(x)M^{c}(\ud s,\ud x),\qquad \forall t\ge 0,f\in D(\A),$$
and has covariance measure $\eta(\ud s,\ud x,\ud y)$.
Let \ $(\omega,s,x)\mapsto G(\omega,s,x)$ \ be a predictable function on $\mathcal{W}^{+}_{0}\times \R^{+}\times E$ such that
\begin{equation}\label{condi:integral Mc}
\p_{\mu}\left[\int_{0}^{t}\ud s\int_{E}2b(x)G^{2}(s,x)X_{s}(\ud x)\right]<+\infty,\qquad \forall t\ge 0.
\end{equation}
Then, following \cite[Section 7.3]{Li}, one can define the stochastic integral of $G$ with respect to the martingale measure $M^{c}(\ud s,\ud x)$,
denoted by
$$\int_{0}^{t}\int_{E}G(s,x)M^{c}(\ud s,\ud x),$$
as the unique square integrable c\`{a}dl\`{a}g martingale with quadratic variation
$$2\int_{0}^{t}\langle b G^{2}(s,\cdot),X_{s}\rangle \ud s.$$ We deduce by \eqref{mean ineq} that condition \eqref{condi:integral Mc} is satisfied in the special case where $G(\omega,s,x)=\varphi(s,x)$ for some $\varphi\in  \mathcal{B}_{b}(\R^{+}\times E)$. To simplify notation, we write in the sequel
$$M^{c}_{t}(\varphi)=\int_{0}^{t}\int_{E}\varphi(s,x)M^{c}(\ud s,\ud x).$$
Now we can define
$$M_{t}(\varphi):=\int_{0}^{t}\int_{E}\varphi(s,x)M(\ud s,\ud x):=M^{d}_{t}(\varphi)+M^{c}_{t}(\varphi)$$
for every $\varphi\in  \mathcal{B}_{b}(\R^{+}\times E)$,
where $M^{d}_{t}(\varphi)$ is the unique purely discontinuous martingale whose jumps are indistinguishable from the process $\langle \varphi(s,\cdot),\Delta X_{s}\rangle 1_{\{\Delta X_{s}\not=0\}}$, and $M^{c}_{t}(\varphi)$ is the unique square integrable c\`{a}dl\`{a}g martingale with quadratic variation $2\int_{0}^{t}\langle b\varphi^{2}(s,\cdot),X_{s}\rangle {\rm d}s$.

\medskip

\begin{proposition}\label{prop6.1}
	Suppose \eqref{newA4} holds.
	For every $f\in \mathcal{B}_{b}(E)$, $t\ge 0$ and $\mu\in\me$,
	\begin{equation}
	\langle f,X_{t}\rangle =\langle \mathfrak{P}_{t}f,X_{0}\rangle+\int_{0}^{t}\int_{E}\mathfrak{P}_{t-s}f(x)M(\ud s,\ud x),\qquad \p_{\mu}\mbox{-a.s.}\label{lem6.1.1}
	\end{equation}
\end{proposition}

\proof We first consider $f\in C_{b}(E)$. Take $q_{0}>c_0$ where $c_{0}$ is the positive constant given in \eqref{mean ineq}. Let $U^{q_{0}}$ and $\mathfrak{U}^{q_{0}}$ be the $q_{0}$-resolvent of $(P_{t})_{t\geq 0}$ and $(\mathfrak{P}_{t})_{t\geq 0}$, respectably.  By taking Laplace transforms of both sides of \eqref{mean equation} we get
\begin{equation}\label{lem6.1.4}
\mathfrak{U}^{q_{0}}f(x)=U^{q_{0}}f(x)+U^{q_{0}}(\gamma-a)\mathfrak{U}^{q_{0}}f(x),\qquad x\in E.
\end{equation}
Recall the concatenation process $\hat{\xi}$ defined in the proof of Lemma \ref{lem3.1}. It is known that $\hat{\xi}$ induces the same topology as $\xi$. Moreover,  $\mathfrak{U}^{q_{0}}f$ is an excessive function with respect to $\hat{\xi}$,  and hence is finely continuous relative to $\hat{\xi}$ (or, equivalently, $\xi$).
Thus,
by condition
\eqref{newA4},
\ $(\gamma-a)\mathfrak{U}^{q_{0}}f\in \mathcal{C}^{\xi}_{b}(E)$.
Equation \eqref{lem6.1.4} implies that $\mathfrak{U}^{q_{0}}f\in D(\A)$ and
\begin{eqnarray}
\A \mathfrak{U}^{q_{0}}f(x)&=&\A U^{q_{0}} f(x)+\A U^{q_{0}}(\gamma-a)\mathfrak{U}^{q_{0}}f(x)\nonumber\\
&=&q_{0} U^{q_{0}} f(x)-f(x)+q_{0} U^{q_{0}}(\gamma-a)\mathfrak{U}^{q_{0}}f(x)-(\gamma-a)\mathfrak{U}^{q_{0}}f(x)\nonumber\\
&=&q_{0} \mathfrak{U}^{q_{0}}f(x)-f(x)-(\gamma-a)\mathfrak{U}^{q_{0}}f(x),\nonumber
\end{eqnarray}
or equivalently,
\begin{equation}
(\A+\gamma-a)\mathfrak{U}^{q_{0}}f(x)=q_{0} \mathfrak{U}^{q_{0}}f(x)-f(x),\qquad x\in E.\nonumber
\end{equation}
Then,  by Theorem \ref{them:martingale}(ii)
$$M_{t}(\mathfrak{U}^{q_{0}}f)=\langle \mathfrak{U}^{q_{0}}f,X_{t}\rangle -\langle \mathfrak{U}^{q_{0}}f,X_{0}\rangle -\int_{0}^{t}\langle q_{0}\mathfrak{U}^{q_{0}}f-f,X_{s}\rangle \ud s$$
is a c\`{a}dl\`{a}g martingale. Using this martingale, one can apply the argument in the proof of \cite[Proposition 2.13]{F1992} with minor modification to show that
\eqref{lem6.1.1} holds for $f\in C_{b}(E)$.

Let $\mathcal{G}$ be the class of bounded measurable functions for which \eqref{lem6.1.1} holds. The above argument shows that $C_{b}(E)\subseteq \mathcal{G}$. By the modified monotone class theorem
(cf.
\cite[Proposition A.2]{Li}), it suffices to prove that $\mathcal{G}$ is closed under bounded pointwise convergence.
Suppose that
$\{f_{n}:n\ge 1\}$ is a sequence of functions from $\mathcal{G}$ and $f$ is the bounded pointwise limit of $f_{n}$. One can easily deduce by bounded convergence theorem that for every $t\ge 0$, \ $\langle f_{n},X_{t}\rangle \to \langle f,X_{t}\rangle$, \ $\langle \mathfrak{P}_{t}f_{n},X_{0}\rangle \to \langle \mathfrak{P}_{t}f,X_{0}\rangle$ \ and that \ $(s,x)\mapsto \1_{\{s\le t\}}\mathfrak{P}_{t-s}f(x)$ \ is the bounded pointwise limit of \ $(s,x)\mapsto \1_{\{s\le t\}}\mathfrak{P}_{t-s}f_{n}(x)$. \ Note that
\begin{eqnarray}
&&\p_{\mu}\left[\left|\int_{0}^{t}\int_{E}\left(\mathfrak{P}_{t-s}f_{n}(x)-\mathfrak{P}_{t-s}f(x)\right)M^{c}(\ud s,\ud x)\right|^{2}\right]\nonumber\\
&=&\p_{\mu}\left[\int_{0}^{t}\int_{E}\langle 2b(\mathfrak{P}_{t-s}f_{n}-\mathfrak{P}_{t-s}f)^{2},X_{s}\rangle \ud s\right]\nonumber\\
&\le&2\|b\|_{\infty}\int_{0}^{t}\Big\langle \mathfrak{P}_{s}\left[\left(\mathfrak{P}_{t-s}f_{n}-\mathfrak{P}_{t-s}f\right)^{2}\right],\mu\Big\rangle \ud s.\nonumber
\end{eqnarray}
By \eqref{mean ineq} and the bounded convergence theorem,
the integral on the right hand side of the above inequality
converges to $0$ as $n\to +\infty$. Hence we get
\begin{equation}\label{lem6.1.2}
\int_{0}^{t}\int_{E}\mathfrak{P}_{t-s}f(x)M^{c}(\ud s,\ud x)=\lim_{n\to+\infty}\int_{0}^{t}\int_{E}\mathfrak{P}_{t-s}f_{n}(x)M^{c}(\ud s,\ud x) \mbox{ in }L^{2}(\p_{\mu}).
\end{equation}
We write $M^{d,n}_{t}$ for \ $\int_{0}^{t}\int_{E}\left(\mathfrak{P}_{t-s}f_{n}(x)-\mathfrak{P}_{t-s}f(x)\right)M^{d}(\ud s,\ud x)$.  This is a purely discontinuous local martingale whose jumps are indistinguishable from the process \ $\langle \mathfrak{P}_{t-s}f_{n}-\mathfrak{P}_{t-s}f,\Delta X_{s}\rangle\1_{\{\Delta X_{s}\not=0\}}$. By
the Burkholder-Davis-Gundy inequality
we have
\begin{eqnarray}
\p_{\mu}\left[|M^{d,n}_{t}|\right]&\le&\p_{\mu}\left[\left(\sum_{s\le t}|\Delta M^{d,n}_{s}|^{2}\right)^{1/2}\right]\nonumber\\
&\le &\p_{\mu}\left[\left(\sum_{s\le t}\Big\langle \left|\mathfrak{P}_{t-s}f_{n}-\mathfrak{P}_{t-s}f\right|,\Delta X_{s}\Big\rangle^{2} \right)^{1/2}\right]\nonumber.
\end{eqnarray}
Applying similar calculations as in \eqref{eq3.15}, we can show that the expectation on the right hand side is less than or equal to
\begin{multline}
\p_{\mu}\left[\int_{0}^{t}\ud s\int_{E}X_{s-}(\ud x)\int_{\me^{0}}\nu\left(\left|\mathfrak{P}_{t-s}f_{n}-\mathfrak{P}_{t-s}f\right|\right)^{2}\1_{\{\nu(1)\le 1\}}H(x,\ud \nu)\right]^{1/2}\nonumber\\
+\p_{\mu}\left[\int_{0}^{t}\ud s\int_{E}X_{s-}(\ud x)\int_{\me^{0}}\nu\left(\left|\mathfrak{P}_{t-s}f_{n}-\mathfrak{P}_{t-s}f\right|\right)\1_{\{\nu(1)>1\}}H(x,\ud \nu)\right].\nonumber
\end{multline}
Note that $\nu(1)\wedge \nu(1)^{2}H(x,\ud \nu)$ is a bounded kernel from $E$ to $\me^{0}$. In view of this and \eqref{mean ineq}, one can show by the bounded convergence theorem that the above two expectations converge to $0$ as $n\to +\infty$. Hence we have
\begin{equation}\label{lem6.1.3}
\int_{0}^{t}\int_{E}\mathfrak{P}_{t-s}f(x)M^{d}(\ud s,\ud x)=\lim_{n\to+\infty}\int_{0}^{t}\int_{E}\mathfrak{P}_{t-s}f_{n}(x)M^{d}(\ud s,\ud x) \mbox{ in }L^{1}(\p_{\mu}).
\end{equation}
\eqref{lem6.1.2} and \eqref{lem6.1.3} imply that there is a subsequence $\{f_{i_{n}}:n\ge 1\}$ such that
\begin{equation}
\lim_{n\to+\infty}\int_{0}^{t}\int_{E}\mathfrak{P}_{t-s}f_{i_{n}}(x)M(\ud s,\ud x)=\int_{0}^{t}\int_{E}\mathfrak{P}_{t-s}f(x)M(\ud s,\ud x)\quad \p_{\mu}\mbox{-a.s.}\label{lem6.1.5}
\end{equation}
Since \eqref{lem6.1.1} holds for $f$ replaced by $f_{i_{n}}$, by letting $n\to +\infty$ we can show by \eqref{lem6.1.5} that it also holds for $f$, and hence $f\in\mathcal{G}$. Therefore $\mathcal{G}$ is closed under bounded pointwise convergence. We complete the proof.\qed

\section{Proofs of the main results}\label{sec:proofs}

\subsection{Interpretation of $\widetilde{P}_{t}$}\label{sec:interpretation}

The following proposition gathers what was already established
in \cite{RSY}. These facts will be used
later in the proofs of the main results.

\begin{proposition}\label{prop1}
Suppose \eqref{newA1}-\eqref{newA3} hold.
For every $x\in E$, define
$$q(x):=\frac{\gamma(x,h)}{h(x)}.$$
Then
$$
{\tt H}_{t}:=\exp\left( \lambda_{1}t-\int_{0}^{t}a(\xi_{s})\ud s
+\int_{0}^{t}q(\xi_{s})\ud s \right)
\frac{h(\xi_{t})}{h(\xi_{0})},\qquad \forall t\ge 0
$$
 is a positive $\Pi_{x}$-martingale with respect to the filtration $\{\mathcal{H}_{t}:t\ge 0\}$. Consequently, the formula
$$
\ud\Pi^{h}_{x}={\tt H}_{t}\,\ud\Pi_{x}\quad
\mbox{on }\mathcal{H}_{t}\cap\{t<\zeta\},\qquad\forall x\in E,
$$
uniquely determines a family of probability measures
$\{\Pi^{h}_{x}:x\in E\}$ on $(\Omega,\mathcal{H})$.
The process $\xi$ under $\{\Pi^{h}_{x}:x\in E\}$ will be
denoted by $\xi^{h}$. The process $\xi^{h}$
is a conservative and recurrent (in the sense of \cite{FOT})
symmetric right Markov process on $E$ with respect to the probability measure \ $\widetilde{m}(\ud y):=h(y)^{2}m(\ud y)$. \
Let $P^{h}_{t}$ denote its transition semigroup, it satisfies that
\begin{eqnarray}
P^{h}_{t}f(x)=\Pi^{h}_{x}\left[f(\xi_{t})\right]
=\frac{e^{\lambda_{1}t}}{h(x)}\Pi_{x}\left[e_{a-q}(t)h(\xi_{t})f(\xi_{t})\right],\nonumber
\end{eqnarray}
for every $x\in E$, $t\ge 0$ and $f\in  \mathcal{B}_{b}(E)$.
Moreover, $\xi^{h}$ has a transition density function with respect to
$\widetilde{m}$.

\end{proposition}

Suppose
$\widehat{\xi}:=((\widehat{\xi}_{t})_{t\ge 0};\widehat{\Pi}^{h}_{x})$
is the $e_{q}(t)$-subprocess of $\xi^{h}$, that is,
$$\widehat{\Pi}^{h}_{x}\left(\widehat{\xi}_{t}\in B\right)=\Pi^{h}_{x}\left[e_{q}(t)\1_{\{\xi_{t}\in B\}}\right],\qquad\forall \ t\ge 0,\ B\in \mathcal{B}(E).$$
In fact, a version of the $e_{q}(t)$-subprocess can be obtained by the following method of curtailment of the lifetime.
Let $Z$ be an exponential random variable of parameter 1 independent of
$\xi^{h}$.
Put
$$\widehat{\zeta}(\omega):=\inf\left\{t\ge 0:\ \int_{0}^{t}q\left(\xi^{h}_{s}(\omega)\right)\ud s\ge Z(\omega)\right\}
\ (=+\infty, \mbox{ if such $t$ does not exist})
,$$
and
\begin{equation} \nonumber
\widehat{\xi}_{t}(\omega):=\begin{cases}
         \xi^{h}_{t}(\omega)
         \quad &\hbox{if } t<\widehat{\zeta}(\omega) ,  \\
         \partial
         \quad &\hbox{if } t\ge \widehat{\zeta}(\omega).
        \end{cases}
\end{equation}
 Then the process $(\widehat{\xi}_{t})_{t\ge 0}$ is equal in law to the $e_{q}(t)$-subprocess of $\xi^{h}$.
Now, we define
\begin{equation}
\kappa(x,\ud y):=	\frac{h(y)\gamma(x,\ud y)}{\gamma(x,h)}\1_{\{\gamma(x,1)>0\}}+\delta_{x}(\ud y)\1_{\{\gamma(x,1)=0\}}, \qquad \mbox{for }x\in E.\label{def:kappa}
\end{equation}
We note that
$\kappa(x,\ud y)$ is a probability kernel on $E$. Let
$\widetilde{\xi}:=((\widetilde{\xi}_{t})_{t\ge 0},\widetilde{\Pi}_{x})$
be the right process constructed from $\widehat{\xi}$ and the instantaneous distribution
\ $\kappa(\widehat{\xi}_{\widehat{\zeta}-}(\omega),\ud y)$ \
by using the so-called ``piecing out" procedure (cf. Ikeda et al.
\cite{INW}), which can be described as follows:
the process $\widetilde{\xi}$ evolves as a copy of $\widehat{\xi}$ until time $\widehat{\zeta}-$, then it is stopped at time $\widehat{\zeta}$ and instantaneously revived by the kernel $\kappa(x,\ud y)$ in the following way:
 At time $\widehat{\zeta}$, the process $\widetilde{\xi}$ is immediately restarted at a new position $y$ which is randomly chosen according to the probability distribution \ $\kappa(\widehat{\xi}_{\widehat{\zeta}-}(\omega),\ud y)$. \
Starting from $y$,\ $\widetilde{\xi}$ evolves again as a copy of $\widehat{\xi}$
and so on, until a countably infinite number of revivals have occurred.
Let $\widetilde{P}_{t}$ be the transition semigroup of $\widetilde{\xi}$.
Naturally by construction it satisfies the renewal equation
\begin{equation}
\widetilde{P}_{t}f(x)=\Pi_{x}^{h}\left[e_{q}(t)f(\xi_{t})\right]+\Pi_{x}^{h}\left[\int_{0}^{t}q(\xi_{s})e_{q}(s)
\kappa(\xi_{s},\widetilde{P}_{t-s}f)\ud s\right]\nonumber
\end{equation}
for every $f\in \mathcal{B}^{+}_{b}(E)$.

\begin{proposition}\label{prop2}
Suppose \eqref{newA1}-\eqref{newA3} hold.
Then, $\widetilde{P}_{t}f(x)$ satisfies \eqref{prop2.0} for every $f\in \mathcal{B}^{+}_{b}(E)$, $t\ge 0$ and $x\in E$.
The probability measure
\begin{equation}\nonumber
\rho(\ud y):=h(y)\widehat{h}(y)m(\ud y)
\end{equation}
is an invariant measure for the semigroup $(\widetilde{P}_{t})_{t\ge 0}$.
Moreover, $\widetilde{\xi}$ has a transition density function $\widetilde{p}(t,x,y)$ with respect to the measure $\rho$
which is given by \eqref{prop2.0 density}.

\end{proposition}

\proof This proposition follows in the same way as \cite[Propositions 4.1]{RSY} with $\gamma(x,\ud y)$ and $\pi^{h}(x,\ud y)$ in the proof of \cite{RSY} replaced by $\gamma(x,\ud y)$ given in \eqref{1.1} and $\kappa(x,\ud y)$ given in \eqref{def:kappa}, respectively. We omit the details here.
The explicit form of  $\widetilde{p}(t,x,y)$  follows from the fact that
\begin{align*}
\int_E \widetilde{p}(t,x,y)f(y)\rho(\ud y)&=\widetilde{P}_{t}f(x)=\frac{e^{\lambda_{1}t}}{h(x)}\mathfrak{P}_{t}(fh)(x)\\
&=\frac{e^{\lambda_{1}t}}{h(x)}\int_E \mathfrak{p}(t,x,y)f(y)h(y)m(\ud y)\\
&=\frac{e^{\lambda_{1}t}}{h(x)}\int_E \mathfrak{p}(t,x,y)f(y)\widehat{h}^{-1}(y)\rho(\ud y)
\end{align*}
for every $x\in E$, $t\ge 0$ and $f\in\mathcal{B}_{b}(E)$.\qed

\begin{remark}\label{remark:many to one}
\rm
Formula \eqref{prop2.0} can be written as
\begin{equation}\nonumber
\frac{\p_{\delta_{x}}\left[\langle f h,X_{t}\rangle\right]}{\p_{\delta_{x}}\left[\langle h,X_{t}\rangle\right]}=\widetilde{\Pi}_{x}\left[f(\widetilde{\xi}_{t})\right],
\qquad\mbox{ for }f\in\mathcal{B}^{+}_{b}(E)\mbox{ and } t\ge 0,
\end{equation}
which enables us to calculate the first moment of the superprocess in
terms of an auxiliary process $\widetilde{\xi}$.
This formula is viewed as an analogue of the ``many-to-one" formula for branching Markov processes.
 In particular, when the branching mechanism is purely local,
 the concatenating procedure described below \eqref{def:kappa} does not occur,
 since $\gamma(x,1)=0$ and $\kappa(x,dy)=\delta_{x}(dy)$ for every $x\in E$.
So in this case,
 the  auxiliary process $\widetilde{\xi}$ runs as a copy of the Doob $h$-transformed process $\xi^{h}$. It holds that
$$\widetilde{P}_{t}f(x)=P^{h}_{t}f(x)=\frac{\e^{\lambda_{1}t}}{h(x)}\Pi_{x}\left[e_{a}(t)h(\xi_{t})f(\xi_{t})\right]=\frac{\e^{\lambda_{1}t}}{h(x)}\mathfrak{P}_{t}(fh)(x),$$
for every $x\in E$, $t\ge 0$ and $f\in\mathcal{B}_{b}(E)$.
\end{remark}

\subsection{Proofs of Proposition \ref{prop:martingale} and Theorems \ref{them7.1}-\ref{them:slln}}

\medskip

\noindent\textit{Proof of Proposition \ref{prop:martingale}: }
This proposition can be proved similarly as \cite[Theorem 3.2]{RSY}. We also give details here for completeness.
By the Markov property of $X$, to show $W^{h}_{t}(X)$ is a martingale, it suffices to prove that
\begin{equation}
\mathfrak{P}_{t}h(x)=e^{-\lambda_{1}t}h(x)\quad\forall x\in E,\ t\ge 0.\label{prop1.1}
\end{equation}
Recall from Proposition \ref{prop1} that $\xi^{h}$ is a conservative process with transition semigroup $P^{h}_{t}$. Let $u(t,x):=\Pi_{x}\left[e_{a-q}(t)h(\xi_{t})\right]$.
Then we have
$$1=P^{h}_{t}1(x)=\frac{e^{\lambda_{1}t}}{h(x)}u(t,x)\quad\forall x\in E,\ t\ge 0,$$
and consequently, $u(t,x)=\mathrm{e}^{-\lambda_{1}t}h(x)$.
Let $A(s,t):=-\int_{s}^{t}(a-q)(\xi_{r}){\rm d}r$.
We note that
$$e^{A(0,t)}-1=-(e^{A(t,t)}-e^{A(0,t)})=\int_{0}^{t}\left(-a(\xi_{s})+q(\xi_{s})\right)e^{A(s,t)}{\rm d}s.$$
Thus by Fubini's theorem and the Markov property of $\xi$, we have
\begin{align*}
u(t,x)&=\Pi_{x}\left[\mathrm{e}^{A(0,t)}h(\xi_{t})\right]\\
&=P_{t}h(x)-\Pi_{x}\left[\int_{0}^{t}a(\xi_{s})e^{A(s,t)}h(\xi_{t}){\rm d}s\right]+\Pi_{x}\left[\int_{0}^{t}q(\xi_{s})e^{A(s,t)}h(\xi_{t}){\rm d}s\right]\\
&=P_{t}h(x)-\Pi_{x}\left[\int_{0}^{t}a(\xi_{s})u(t-s,\xi_{s}){\rm d}s\right]+\Pi_{x}\left[\int_{0}^{t}\frac{\gamma(\xi_{s},h)}{h(\xi_{s})}
u(t-s,\xi_{s}){\rm d}s\right]\nonumber\\
&=P_{t}h(x)-\Pi_{x}\left[\int_{0}^{t}a(\xi_{s})u(t-s,\xi_{s}){\rm d}s\right]+\Pi_{x}\left[\int_{0}^{t}\gamma(\xi_{s},u^{t-s}){\rm d}s\right].\nonumber
\end{align*}
In the last equality we use the fact that $u(t-s,x)=e^{-\lambda_{1}(t-s)}h(x)$ twice.
The above equality implies that $u(t,x)=\mathfrak{P}_{t}h(x)$ is the unique locally bounded solution to \eqref{mean equation} for $f=h$. Hence we prove \eqref{prop1.1}.\qed

\medskip

For the remainder of this section
we assume that conditions \eqref{newA1}-\eqref{newA4} hold and \eqref{newA5} is satisfied for some constant $p\in (1,2]$.
Conditions used in each lemma are stated explicitly.
Let us
explain shortly how to prove Theorems \ref{them7.1}-\ref{them:slln}.

\noindent(i) Since $W_t^h(X)$ is  a martingale, in order to prove Theorem \ref{them7.1},
we shall prove that
	\begin{equation}
	\p_{\mu}\left[\sup_{t\ge 0}W^{h}_{t}(X)^{p}\right]<+\infty, \qquad \forall\mu\in\me.\label{them7.1.1}
	\end{equation}

\noindent(ii) Note that for
any $t,s\ge 0$ and $f\in \mathcal{B}^{+}_{b}(E)$,
\begin{eqnarray}
\e^{\lambda_{1}(t+s)}\langle f,X_{t+s}\rangle&=&\left(\e^{\lambda_{1}(t+s)}\langle f,X_{t+s}\rangle -\p_{\mu}\left[\e^{\lambda_{1}(t+s)}\langle f,X_{t+s}\rangle |\mathcal{F}_{t}\right]\right)\nonumber\\
&&+\p_{\mu}\left[\e^{\lambda_{1}(t+s)}\langle f,X_{t+s}\rangle |\mathcal{F}_{t}\right].\label{eq:4.4}
\end{eqnarray}
We shall prove Theorem \ref{them:wlln} by showing the
$L^{p}$-convergence
of the two summands in \eqref{eq:4.4}.
This is done through
Lemmas \ref{lem7.2} and \ref{lem7.3}.

\noindent(iii) The proof of  Theorem \ref{them:slln} follows two main steps. Firstly we shall prove the almost sure convergence along lattice times (Lemma \ref{lem:wlln}). Then we extend it to continuous time. The transition from discrete to continuous time is obtained through approximation of bounded functions by resolvent functions (Lemma \ref{lem7.5} and equation \eqref{them7.3.2}).

\medskip

\noindent\textit{Proof of Theorem \ref{them7.1}}:
Fix an arbitrary $\mu\in\me$.
By Propositions
\ref{prop:martingale} and \ref{prop6.1}
we have
$$W^{h}_{t}(X)=\langle h,X_{0}\rangle +\int_{0}^{t}\int_{E}\e^{\lambda_{1}s}h(x)M(\ud s,\ud x).$$
By Doob's martingale inequality and Jensen's inequality we have
\begin{eqnarray}
&&\p_{\mu}\left[\sup_{0\le r\le t}\left|\int_{0}^{r}\int_{E}\e^{\lambda_{1}s}h(x)M^{c}(\ud s,\ud x)\right|^{p}\right]^{2/p}\nonumber\\
&\le&\left(\frac{p}{p-1}\right)^{2}\p_{\mu}\left[\left|\int_{0}^{t}\int_{E}\e^{\lambda_{1}s}h(x)M^{c}(\ud s,\ud x)\right|^{p}\right]^{2/p}\nonumber\\
&\le&\left(\frac{p}{p-1}\right)^{2}\p_{\mu}\left[\left|\int_{0}^{t}\int_{E}\e^{\lambda_{1}s}h(x)M^{c}(\ud s,\ud x)\right|^{2}\right]\nonumber\\
&=&2\left(\frac{p}{p-1}\right)^{2}\p_{\mu}\left[\int_{0}^{t}\e^{2\lambda_{1}s}\langle bh^{2},X_{s}\rangle \ud s\right]\nonumber\\
&\le&2\left(\frac{p}{p-1}\right)^{2}\|b\|_{\infty}\|h\|_{\infty}\langle h,\mu\rangle \int_{0}^{t}\e^{\lambda_{1}s} \ud s.\nonumber
\end{eqnarray}
Since $\lambda_{1}<0$, by letting $t\to +\infty$, we get
\begin{equation}
\p_{\mu}\left[\sup_{r\ge 0}\left|\int_{0}^{r}\int_{E}\e^{\lambda_{1}s}h(x)M^{c}(\ud s,\ud x)\right|^{p}\right]^{2/p}<+\infty.\label{them7.1.2}
\end{equation}
We note that \ $t\mapsto \int_{0}^{t}\int_{E}\e^{\lambda_{1}s}h(x)M^{d}(\ud s,\ud x)$ \ is a purely discontinuous local martingale whose jumps are
indistinguishable
from the process $\e^{\lambda_{1}s}\langle h,\Delta X_{s}\rangle \1_{\{\Delta X_{s}\not=0\}}$.
Hence, by
the Doob's martingale inequality and
the Burkholder-Davis-Gundy inequality
for purely discontinuous local martingale, we have
\begin{eqnarray}
&&\p_{\mu}\left[\sup_{0\le r\le t}\left|\int_{0}^{r}\int_{E}\e^{\lambda_{1}s}h(x)M^{d}(\ud s,\ud x)\right|^{p}\right]\nonumber\\
&\le&\left(\frac{p}{p-1}\right)^{p}\p_{\mu}\left[\left|\int_{0}^{t}\int_{E}\e^{\lambda_{1}s}h(x)M^{d}(\ud s,\ud x)\right|^{p}\right]\nonumber
\end{eqnarray}
\begin{eqnarray}
&\le& c_{1}\p_{\mu}\left[\left(\sum_{0\le s\le t}\e^{2\lambda_{1}s}\langle h,\Delta X_{s}\rangle^{2}\right)^{p/2}\right]\nonumber\\
&\le& c_{1}\p_{\mu}\left[\sum_{0\le s\le t}\e^{\lambda_{1}p s}\langle h,\Delta X_{s}\rangle^{p}\right]\nonumber\\
&=&c_{1}\p_{\mu}\left[\int_{0}^{t}\e^{\lambda_{1}ps}\ud s\int_{E}X_{s-}(\ud x)\int_{\me^{0}}\nu(h)^{p}H(x,\ud \nu)\right],\nonumber
\end{eqnarray}
where $c_{1}=c_{1}(p)$ is a positive constant. In the third inequality we use the fact that $(a_{1}+\cdots+a_{n})^{p/2}\le a_{1}^{p/2}+\cdots+a_{n}^{p/2}$ for any $a_{1},\ldots,a_{n}\ge 0$ and $p\in(1,2]$.
By Proposition \ref{prop2},
the expectation on the right hand side
is equal to
\begin{eqnarray}
&&\int_{0}^{t}\e^{\lambda_{1}ps}\Bigg\langle \mathfrak{P}_{s}\left(\int_{\me^{0}}\nu(h)^{p}H(\cdot,\ud \nu)\right),\mu\Bigg\rangle \ud s\nonumber\\
&=&\int_{0}^{t}\e^{\lambda_{1}(p-1)s}\Bigg\langle h \widetilde{P}_{s}\left(h^{-1}\int_{\me^{0}}\nu(h)^{p}H(\cdot,\ud \nu)\right),\mu \Bigg\rangle \ud s\nonumber\\
&\le&\left\|h^{-1}\int_{\me^{0}}\nu(h)^{p}H(\cdot,\ud \nu)\right\|_{\infty}\langle h,\mu\rangle\int_{0}^{t}\e^{\lambda_{1}(p-1)s}\ud s.\nonumber
\end{eqnarray}
Letting $t\to+\infty$, we get
\begin{equation}
\p_{\mu}\left[\sup_{s\ge 0}\left|\int_{0}^{s}\int_{E}\e^{\lambda_{1}s}h(x)M^{d}(\ud s,\ud x)\right|^{p}\right]<+\infty.\label{them7.1.3}
\end{equation}
Therefore,
\eqref{them7.1.1} follows directly from \eqref{them7.1.2} and \eqref{them7.1.3}. We complete the proof.\qed

\medskip

In order to
simplified computations,
we will work with
the test functions
$f=\phi h$ with $ \phi\in \mathcal{B}^{+}_{b}(E)$.

\begin{lemma}\label{lem7.2}
Suppose \eqref{newA1}-\eqref{newA6} hold.
For any $\phi\in \mathcal{B}^{+}_{b}(E)$ and $\mu\in\me$,
\begin{equation}\nonumber
\lim_{s\to+\infty}\lim_{t\to +\infty}\p_{\mu}\left[\e^{\lambda_{1}(t+s)}\langle \phi h,X_{t+s}\rangle|\mathcal{F}_{t}\right]=(\phi h,\hat{h}) W^{h}_{\infty}(X)\qquad \p_{\mu}\mbox{-a.s. and in }L^{p}(\p_{\mu}).
\end{equation}
\end{lemma}

\proof Fix $\phi\in \mathcal{B}^{+}_{b}(E)$ and $\mu\in\me$. By Theorem \ref{them7.1}, it suffices to prove that
\begin{equation}\label{lem7.2.1}
\lim_{s\to+\infty}\lim_{t\to +\infty}\left(\p_{\mu}\left[\e^{\lambda_{1}(t+s)}\langle \phi h,X_{t+s}\rangle|\mathcal{F}_{t}\right]-(\phi h,\hat{h}) W^{h}_{t}(X)\right)=0\ \p_{\mu}\mbox{-a.s. and in }L^{p}(\p_{\mu}).
\end{equation}
For any $s,t\ge 0$, by the Markov property
\begin{eqnarray}
\p_{\mu}\left[\e^{\lambda_{1}(t+s)}\langle \phi h,X_{t+s}\rangle|\mathcal{F}_{t}\right]
&=&\e^{\lambda_{1}(t+s)}\p_{X_{t}}\left[\langle \phi h,X_{s}\rangle\right]\nonumber\\
&=&\e^{\lambda_{1}(t+s)}\langle \mathfrak{P}_{s}(\phi h),X_{t}\rangle\nonumber\\
&=&\e^{\lambda_{1}t}\langle h\widetilde{P}_{s}\phi,X_{t}\rangle.\nonumber
\end{eqnarray}
Hence,
$$\left|\p_{\mu}\left[\e^{\lambda_{1}(t+s)}\langle \phi h,X_{t+s}\rangle|\mathcal{F}_{t}\right]-(\phi h,\hat{h}) W^{h}_{t}(X)\right|
\le \e^{\lambda_{1}t}\Big\langle h\left|\widetilde{P}_{s}\phi-\langle \phi h,\hat{h}\rangle\right|,X_{t}\Big\rangle.$$
It follows by condition \eqref{newA6} and Proposition \ref{prop2} that for any $\epsilon>0$ any $s$ sufficiently large,
$$\sup_{x\in E}\left|\widetilde{P}_{s}\phi(x)-\langle \phi h,\hat{h}\rangle\right|\le \sup_{x\in E}\int_{E}|\phi(y)| \left|\widetilde{p}(s,x,y)-1\right|\rho(\ud y)\le \epsilon\|\phi\|_{\infty},$$
in which case
\begin{equation}\label{lem7.2.2}
\left|\p_{\mu}\left[\e^{\lambda_{1}(t+s)}\langle \phi h,X_{t+s}\rangle|\mathcal{F}_{t}\right]-(\phi h,\hat{h}) W^{h}_{t}(X)\right|\le \epsilon \|\phi\|_{\infty}W^{h}_{t}(X)
\end{equation}
Since $W^{h}_{t}(X)\to W^{h}_{\infty}(X)$ $\p_{\mu}$-a.s. and in $L^{p}(\p_{\mu})$, we get \eqref{lem7.2.1} by letting $t\to+\infty$ and $\epsilon\to 0$ in \eqref{lem7.2.2}.\qed

\begin{lemma}\label{lem7.3}
Suppose \eqref{newA1}-\eqref{newA6} hold.
For any $\phi\in \mathcal{B}^{+}_{b}(E)$, $\mu\in \me$
and $s\ge 0$,
\begin{equation*}
\lim_{t\to +\infty}\left(\e^{\lambda_{1}(t+s)}\langle \phi h,X_{t+s}\rangle-\p_{\mu}\left[\e^{\lambda_{1}(t+s)}\langle \phi h,X_{t+s}\rangle|\mathcal{F}_{t}\right]\right)=0\quad \mbox{ in }L^{p}(\p_{\mu}).
\end{equation*}
Moreover, for any $m\in\mathbb{N}$ and $\sigma>0$ the following holds $\p_{\mu}$-a.s. and in $L^{p}(\p_{\mu})$.
\begin{equation*}
\lim_{n\to+\infty}\left(\e^{\lambda_{1}(m+n)\sigma}\langle \phi h,X_{(m+n)\sigma}\rangle-\p_{\mu}\left[\e^{\lambda_{1}(m+n)\sigma}\langle \phi h,X_{(m+n)\sigma}\rangle|\mathcal{F}_{n\sigma}\right]\right)=0.
\end{equation*}
\end{lemma}

\proof For any $T>0$, we define
\begin{equation*}
L^{T}_{s,t}(\phi):=\int_{s}^{t}\int_{E}\mathfrak{P}_{T-r}(\phi h)(x)M(\ud r,\ud x),
\qquad 0\le s\le t\le T, \ \phi\in \mathcal{B}^{+}_{b}(E),
\end{equation*}
and define \ $L^{T,c}_{s,t}(\phi)$ \ and \ $L^{T,d}_{s,t}(\phi)$ \ analogously with $M$ replaced by $M^{c}$ and $M^{d}$ respectively.
For simplicity, $T$ is omitted when $t=T$.
By the fact that $(a+b)^{p}\le 2^{p-1}(a^{p}+b^{p})$
for any $a,b\ge 0$, $p\in (1,2]$,
and Jensen's inequality, we have
$$\p_{\mu}\left[|L^{T}_{0,t}(\phi)|^{p}\right]\le 2^{p-1}\left[\p_{\mu}\left[L^{T,c}_{0,t}(\phi)^{2}\right]^{p/2}+\p_{\mu}\left[|L^{T,d}_{0,t}(\phi)|^{p}\right]\right].$$
Applying similar calculations as in the proof of Theorem \ref{them7.1}, one can show that the two expectations
on the right hand side of the above inequality
are finite and hence $\p_{\mu}\left[|L^{T}_{0,t}(\phi)|^{p}\right]<+\infty$ for every $t\in [0,T]$.
Thus the local martingale \ $[0,T]\ni t\mapsto L^{T}_{0,t}(\phi)$ \ is an $L^{p}$-integrable martingale. Using this and Proposition \ref{prop6.1}, we have
\begin{eqnarray}
&&\e^{\lambda_{1}(t+s)}\langle \phi h,X_{t+s}\rangle-\p_{\mu}\left[\e^{\lambda_{1}(t+s)}\langle \phi h,X_{t+s}\rangle|\mathcal{F}_{t}\right]\nonumber\\
&=&\e^{\lambda_{1}(t+s)}L_{0,t+s}(\phi)-\e^{\lambda_{1}(t+s)}\p_{\mu}\left[L_{0,t+s}(\phi)|\mathcal{F}_{t}\right]\nonumber\\
&=&\e^{\lambda_{1}(t+s)}L_{0,t+s}(\phi)-\e^{\lambda_{1}(t+s)}L^{t+s}_{0,t}(\phi)\nonumber\\
&=&\e^{\lambda_{1}(t+s)}L_{t,t+s}(\phi).\label{lem7.3.1}
\end{eqnarray}
On one hand, by Proposition \ref{prop2}, we have
\begin{eqnarray}
\p_{\mu}\left[\left(\e^{\lambda_{1}(t+s)}L^{c}_{t,t+s}(\phi)\right)^{2}\right]
&=&2\p_{\mu}\left[\int_{t}^{t+s}\e^{2\lambda_{1}(t+s)}\langle b\mathfrak{P}_{t+s-r}(\phi h)^{2},X_{r}\rangle \ud r\right]\nonumber\\
&=&2\p_{\mu}\left[\int_{t}^{t+s}\e^{2\lambda_{1}r}\langle bh^{2}(\widetilde{P}_{t+s-r}\phi)^{2},X_{r}\rangle \ud r\right]\nonumber\\
&\le&2\|b\|_{\infty}\|h\|_{\infty}\|\phi\|_{\infty}\p_{\mu}\left[\int_{t}^{t+s}\e^{2\lambda_{1}r}\langle h\widetilde{P}_{t+s-r}\phi,X_{r}\rangle \ud r\right]\nonumber\\
&=&2\|b\|_{\infty}\|h\|_{\infty}\|\phi\|_{\infty}\langle h\widetilde{P}_{t+s}\phi,\mu\rangle\int_{t}^{t+s}\e^{\lambda_{1}r}\ud r.\nonumber
\end{eqnarray}
Immediately,
$$\p_{\mu}\left[\left(\e^{\lambda_{1}(t+s)}L^{c}_{t,t+s}(\phi)\right)^{2}\right]\le 2\|b\|_{\infty}\|\phi\|^{2}_{\infty}\|h\|_{\infty}\langle h,\mu\rangle \int_{t}^{t+s}\e^{\lambda_{1}r}\ud r\quad \forall t,s\ge 0.$$
Since $\lambda_{1}<0$, it follows that
$$\lim_{t\to +\infty}\p_{\mu}\left[\left(\e^{\lambda_{1}(t+s)}L^{c}_{t,t+s}(\phi)\right)^{2}\right]=0,\ \sum_{n=1}^{+\infty}\p_{\mu}\left[\left(\e^{\lambda_{1}(m+n)\sigma}L^{c}_{n\sigma,(m+n)\sigma}(\phi)\right)^{2}\right]<+\infty,$$
which in turn implies that
\begin{eqnarray}
&&\lim_{t\to +\infty}\e^{\lambda_{1}(t+s)}L^{c}_{t,t+s}(\phi)=0\mbox{ in }L^{2}(\p_{\mu}),\nonumber\\ &&\lim_{n\to+\infty}\e^{\lambda_{1}(m+n)\sigma}L^{c}_{n\sigma,(m+n)\sigma}(\phi)=0\quad
\p_{\mu}\mbox{-a.s. and in }L^{2}(\p_{\mu}).
\label{lem7.3.2}
\end{eqnarray}
On the other hand, by
the Burkholder-Davis-Gundy inequality
and the fact that $(a+b)^{p/2}\le a^{p/2}+b^{p/2}$ for any $a,b\ge 0$ and $p\in (1,2]$, we have that
\begin{align*}
&\hspace{0.5cm}\p_{\mu}\left[\left|\e^{\lambda_{1}(t+s)}L^{d}_{t,t+s}(\phi)\right|^{p}\right]\nonumber\\
&\le c_{1}\p_{\mu}\left[\left(\sum_{t\le r\le t+s}\e^{2\lambda_{1}(t+s)}\langle \mathfrak{P}_{t+s-r}(\phi h),\Delta X_{r}\rangle^{2}\right)^{p/2}\right]\nonumber\\
&\le c_{1}\p_{\mu}\left[\sum_{t\le r\le t+s}\e^{\lambda_{1}p(t+s)}\langle \mathfrak{P}_{t+s-r}(\phi h),\Delta X_{r}\rangle^{p}\right]\nonumber\\
&= c_{1}\p_{\mu}\left[\int_{t}^{t+s}\e^{\lambda_{1}p(t+s)}\ud r\int_{E}X_{r-}(\ud x)\int_{\me^{0}}\nu(\mathfrak{P}_{t+s-r}(\phi h))^{p}H(x,\ud \nu)\right]\nonumber\\
&= c_{1}\int_{t}^{t+s}\e^{\lambda_{1}(p-1)r}\Bigg\langle h\widetilde{P}_{r}\left(h^{-1}\int_{\me^{0}}\nu\left(h\widetilde{P}_{t+s-r}\phi\right)^{p}H(\cdot,\ud \nu)\right),\mu\Bigg\rangle\ \ud r,
\end{align*}
where $c_{1}=c_{1}(p)$ is a positive constant. Since $\|\widetilde{P}_{t+s-r}\phi\|_{\infty}\le \|\phi\|_{\infty}$, we get
\begin{equation}
\p_{\mu}\left[\left|\e^{\lambda_{1}(t+s)}L^{d}_{t,t+s}(\phi)\right|^{p}\right]\le c_{1}\|\phi\|^{p}_{\infty}\int_{t}^{t+s}\e^{\lambda_{1}(p-1) r}\langle h\widetilde{P}_{r}F,\mu\rangle \ud r\label{lem7.3.3}
\end{equation}
where \ $F(x):=h^{-1}(x)\int_{\me^{0}}\nu(h)^{p}H(x,\ud \nu)$. \ It follows by condition \eqref{newA6} that for any $\epsilon>0$ there exists $t$ sufficiently large such that
$$\sup_{x\in E}\left|\widetilde{P}_{r}F(x)\right|\le (1+\epsilon)\rho(F),\qquad \forall r\ge t,$$
in which case the integral
on the right hand side
of \eqref{lem7.3.3} is
less than or equal to
$$\langle h,\mu\rangle (1+\epsilon)\rho(F)\int_{t}^{t+s}\e^{\lambda_{1}(p-1)r}\ud r.$$
Thus, we get
$$\lim_{t\to+\infty}\p_{\mu}\left[\left|\e^{\lambda_{1}(t+s)}L^{d}_{t,t+s}(\phi)\right|^{p}\right]=0,\  \sum_{n=1}^{+\infty}\p_{\mu}\left[\left|\e^{\lambda_{1}(m+n)\sigma}L^{d}_{n\sigma,(n+m)\sigma}(\phi)\right|^{p}\right]<+\infty,$$
and consequently, by Borel-Cantelli lemma
\begin{eqnarray}
&&\lim_{t\to+\infty}\e^{\lambda_{1}(t+s)}L^{d}_{t,t+s}(\phi)=0\mbox{ in }L^{p}(\p_{\mu}),\nonumber\\ &&\lim_{n\to+\infty}\e^{\lambda_{1}(m+n)\sigma}L^{d}_{n\sigma,(n+m)\sigma}(\phi)=0\quad
\p_{\mu}\mbox{-a.s. and in }L^{p}(\p_{\mu}).
\label{lem7.3.4}
\end{eqnarray}
The lemma follows from \eqref{lem7.3.1}, \eqref{lem7.3.2} and \eqref{lem7.3.4}.\qed

\medskip

\noindent\textit{Proof of Theorem \ref{them:wlln}: } In view of \eqref{eq:4.4}, this theorem is an immediate consequence of Lemma \ref{lem7.2} and Lemma \ref{lem7.3}.\qed

\medskip

Lemma \ref{lem7.2} and Lemma \ref{lem7.3}
also give
the following result.

\begin{lemma}\label{lem:wlln}
Suppose \eqref{newA1}-\eqref{newA6} hold. Then, for any $\sigma>0$, $\mu\in\me$ and $f\in \mathcal{B}^{+}(E)$ with $f/h$ bounded,
$$\lim_{n\to+\infty}\e^{\lambda_{1}n\sigma}\langle f,X_{n\sigma}\rangle =(f,\hat{h}) W^{h}_{\infty}(X)\qquad \p_{\mu}\mbox{-a.s.}$$
\end{lemma}

This result is the SLLN for superprocesses along lattice times. For the transition to continuous time, we need the following lemma.

\medskip

\begin{lemma}\label{lem7.5}
Suppose \eqref{newA1}-\eqref{newA6} hold.
For any $\phi \in \mathcal{B}^{+}_{b}(E)$, $\alpha >0$ and $\mu\in\me$,
\begin{equation}\nonumber
\lim_{t\to+\infty}\e^{\lambda_{1}t}\langle (\alpha \widetilde{U}^{\alpha }\phi)h,X_{t}\rangle =(\phi h,\hat{h}) W^{h}_{\infty}(X)\qquad \p_{\mu}\mbox{-a.s.}
\end{equation}
where \ $\widetilde{U}^{\alpha }\phi(x):=\int_{0}^{+\infty}\e^{-\alpha  s}\widetilde{P}_{s}\phi(x)\ud s$.
\end{lemma}

\proof Fix $\phi \in \mathcal{B}^{+}_{b}(E)$ and $\alpha ,\sigma>0$. Let $g(x):=\alpha \widetilde{U}^{\alpha }\phi(x)$ for $x\in E$. Suppose $t\in [n\sigma,(n+1)\sigma)$. We have
\begin{eqnarray}
&&\e^{\lambda_{1}t}\langle gh,X_{t}\rangle-(\phi h,\hat{h}) W^{h}_{\infty}(X)\nonumber\\
&=&\left(\e^{\lambda_{1}t}\langle gh,X_{t}\rangle-\p_{\mu}\left[\e^{\lambda_{1}(n+1)\sigma}\langle gh,X_{(n+1)\sigma}\rangle|\mathcal{F}_{t}\right]\right)\nonumber\\
&&+\left(\p_{\mu}\left[\e^{\lambda_{1}(n+1)\sigma}\langle gh,X_{(n+1)\sigma}\rangle|\mathcal{F}_{t}\right]-\p_{\mu}\left[\e^{\lambda_{1}(n+1)\sigma}\langle gh,X_{(n+1)\sigma}\rangle|\mathcal{F}_{n\sigma}\right]\right)\nonumber\\
&&+\left(\p_{\mu}\left[\e^{\lambda_{1}(n+1)\sigma}\langle gh,X_{(n+1)\sigma}\rangle|\mathcal{F}_{n\sigma}\right]
-(\phi h,\hat{h}) W^{h}_{\infty}(X)\right)\nonumber\\
&=:&I^{(1)}(t,(n+1)\sigma)+I^{(2)}(t,n\sigma,(n+1)\sigma)+I^{(3)}(n\sigma,(n+1)\sigma).\label{lem7.5.1}
\end{eqnarray}
By Markov property we have
\begin{eqnarray}
I^{(1)}(t,(n+1)\sigma)&=&\e^{\lambda_{1}t}\langle gh,X_{t}\rangle-\e^{\lambda_{1}(n+1)\sigma}\p_{X_{t}}\left[\langle gh,X_{(n+1)\sigma-t}\rangle\right]\nonumber\\
&=&\e^{\lambda_{1}t}\langle gh,X_{t}\rangle-\e^{\lambda_{1}t}\langle h\widetilde{P}_{(n+1)\sigma-t}(g),X_{t}\rangle\nonumber\\
&=&\e^{\lambda_{1}t}\Big\langle h\left(g-\widetilde{P}_{(n+1)\sigma-t}(g)\right),X_{t}\Big\rangle.\label{lem7.5.2}
\end{eqnarray}
Note that
\begin{align}
g(x)-\widetilde{P}_{(n+1)\sigma-t}&g(x)
=\int_{0}^{+\infty}\hspace{-.5cm}\alpha \e^{-\alpha  s}\widetilde{P}_{s}\phi(x)\ud s-\e^{\alpha ((n+1)\sigma-t)}\int_{(n+1)\sigma-t}^{+\infty}\alpha \e^{-\alpha  s}\widetilde{P}_{s}\phi(x)\ud s\nonumber\\
&=\left(1-\e^{\alpha ((n+1)\sigma-t)}\right)g(x)+\e^{\alpha ((n+1)\sigma-t)}\int_{0}^{(n+1)\sigma-t}\alpha \e^{-qs}\widetilde{P}_{s}\phi(x)\ud s.\nonumber
\end{align}
Hence $\|g-\widetilde{P}_{(n+1)\sigma-t}g\|_{\infty}\le 2\|\phi \|_{\infty}\left(\e^{\alpha ((n+1)\sigma-t)}-1\right)$, and by \eqref{lem7.5.2}
$$\sup_{t\in [n\sigma,(n+1)\sigma)}|I^{(1)}(t,(n+1)\sigma)|\le 2\|\phi \|_{\infty}\left(\e^{\alpha \sigma}-1\right)\sup_{t\in [n\sigma,(n+1)\sigma)}W^{h}_{t}(X).$$
Since $W_{t}^{h}(X)$ converges $\p_{\mu}$-a.s. to a finite limit, the above inequality implies that
\begin{equation}\label{lem7.5.3}
\lim_{\sigma\to 0}\lim_{n\to+\infty}\sup_{t\in [n\sigma,(n+1)\sigma)}|I^{(1)}(t,(n+1)\sigma)|=0\quad \p_{\mu}\mbox{-a.s.}
\end{equation}
Recall the definition of $L^{T}_{s,t}(\phi )$ given in the proof of Lemma \ref{lem7.3}. By Proposition \ref{prop6.1} we have
\begin{align}
I^{(2)}(t,n\sigma,(n+1)\sigma)
=&\ \p_{\mu}\left[\e^{\lambda_{1}(n+1)\sigma}L_{0,(n+1)\sigma}(g)|\mathcal{F}_{t}\right]-
\p_{\mu}\left[\e^{\lambda_{1}(n+1)\sigma}L_{0,(n+1)\sigma}(g)|\mathcal{F}_{n\sigma}\right]\nonumber\\
=&\ \e^{\lambda_{1}(n+1)\sigma}L^{(n+1)\sigma}_{0,t}(g)-\e^{\lambda_{1}(n+1)\sigma}L^{(n+1)\sigma}_{0,n\sigma}(g)\nonumber\\
=&\ \e^{\lambda_{1}(n+1)\sigma}L^{(n+1)\sigma}_{n\sigma,t}(g),\label{lem7.5.4}
\end{align}
where for
the second equality
we use the fact that \ $[0,(n+1)\sigma]\ni t\mapsto L^{(n+1)\sigma}_{0,t}(g)$ \ is a martingale. It
follows from
Doob's martingale inequality
that, for $p\in (1,2]$,
\begin{eqnarray}
&&\p_{\mu}\left[\e^{\lambda_{1}p(n+1)\sigma}\sup_{t\in [n\sigma,(n+1)\sigma)}|L^{(n+1)\sigma}_{n\sigma,t}(g)|^{p}\right]\nonumber\\
&\le&c_{1}\e^{\lambda_{1}p(n+1)\sigma}\p_{\mu}\left[|L_{n\sigma,(n+1)\sigma}(g)|^{p}\right]\nonumber\\
&\le&c_{2}\e^{\lambda_{1}p(n+1)\sigma}\left(\p_{\mu}\left[|L^{c}_{n\sigma,(n+1)\sigma}(g)|^{2}\right]^{p/2}
+\p_{\mu}\left[|L^{d}_{n\sigma,(n+1)\sigma}(g)|^{p}\right]\right),\label{lem7.5.5}
\end{eqnarray}
where $c_{i}=c_{i}(p)>0,\ i=1,2.$ We have showed in the proof of Lemma \ref{lem7.3} that for $n$ sufficiently large,
$$\e^{2\lambda_{1}(n+1)\sigma}\p_{\mu}\left[L^{c}_{n\sigma,(n+1)\sigma}(g)^{2}\right]\le c_{3}\int_{n\sigma}^{(n+1)\sigma}\e^{\lambda_{1}r}\ud r,$$
and
$$\e^{\lambda_{1}p(n+1)\sigma}\p_{\mu}\left[\left|L^{d}_{n\sigma,(n+1)\sigma}(g)\right|^{p}\right]\le c_{4}\int_{n\sigma}^{(n+1)\sigma}\e^{\lambda_{1}(p-1)r}\ud r$$
for some positive constants $c_{i}$, $i=3,4$ independent of $n$.
Using the above estimates and \eqref{lem7.5.5}
one can easily show that
$\sum_{n=N}^{+\infty}\p_{\mu}\left[\e^{\lambda_{1}(n+1)\sigma}\sup_{t\in [n\sigma,(n+1)\sigma)}|L^{(n+1)\sigma}_{n\sigma,t}(g)|^{p}\right]$ is finite for $N$ large enough. Thus by Borel-Cantelli lemma and \eqref{lem7.5.4} we get
\begin{equation}\label{lem7.5.6}
\sup_{t\in [n\sigma,(n+1)\sigma)}\left|I^{(2)}(t,n\sigma,(n+1)\sigma)\right|=\sup_{t\in [n\sigma,(n+1)\sigma)}\e^{\lambda_{1}(n+1)\sigma}|L^{(n+1)\sigma}_{t,n\sigma}(g)|\to 0\quad\p_{\mu}\mbox{-a.s.}
\end{equation}
as $n\to+\infty$. Finally we have
\begin{eqnarray}
I^{(3)}(n\sigma,(n+1)\sigma)&=&\e^{\lambda_{1}(n+1)\sigma}\p_{X_{n\sigma}}\left[\langle gh,X_{\sigma}\rangle\right]-(\phi h,\hat{h}) W^{h}_{\infty}(X)\nonumber\\
&=&\e^{\lambda_{1}n\sigma}\langle h\widetilde{P}_{\sigma}g,X_{n\sigma}\rangle-(\phi h,\hat{h}) W^{h}_{\infty}(X).\nonumber
\end{eqnarray}
Recall that $(\widetilde{P}_{t})_{t\ge 0}$ is invariant with respect to the measure $\rho(\ud x)=h(x)\hat{h}(x)m(\ud x)$. We have $(h\widetilde{P}_{\sigma}g,\hat{h})=\alpha \int_{0}^{\infty}e^{-\alpha  s}(h\widetilde{P}_{\sigma+s}\phi ,\hat{h})\ud s=(\phi h,\hat{h})$.
It then follows by
Lemma \ref{lem:wlln}
that
\begin{equation}
\lim_{n\to+\infty}I^{(3)}(n\sigma,(n+1)\sigma)=0\quad\p_{\mu}\mbox{-a.s.}\label{lem7.5.7}
\end{equation}
In view of \eqref{lem7.5.3}, \eqref{lem7.5.6} and \eqref{lem7.5.7}, one can prove this lemma
by letting first $n\to +\infty$
and then $\sigma\to 0$ in \eqref{lem7.5.1}.\qed

\medskip

\noindent Finally we shall prove SLLN along continuous time under assumptions \eqref{newA1}-\eqref{newA7}.

\medskip

\noindent\textit{Proof of Theorem \ref{them:slln}: }
By \cite[Lemma 7.1]{CRY} it suffices to prove that for any $\mu\in\me$ and $\phi\in C^{+}_{0}(E)$,
\begin{equation}\label{them7.3.1}
\lim_{t\to+\infty}\e^{\lambda_{1}t}\langle \phi h,X_{t}\rangle =(\phi h,\hat{h}) W^{h}_{\infty}(X)\quad \p_{\mu}\mbox{-a.s.}
\end{equation}
Condition \eqref{newA7} implies that for any $\phi\in C^{+}_{0}(E)$
\begin{equation}\label{them7.3.2}
\|\alpha \widetilde{U}^{\alpha}\phi -\phi \|_{\infty}\le \int_{0}^{{\infty}}\alpha e^{-\alpha t }\|\widetilde{P}_{t}\phi -\phi \|_{\infty}{\rm d}t
=\int_{0}^{{\infty}} e^{- s }\|\widetilde{P}_{s/\alpha}\phi -\phi \|_{\infty}{\rm d}s\to 0
\end{equation}
as $\alpha \to{\infty}$. For any $\alpha>0$, we have
\begin{eqnarray}
&&\left|e^{\lambda_{1}t}\langle \phi h,X_{t}\rangle-(\phi h,\hat{h}) W^{h}_{\infty}(X)\right|\nonumber\\
&\le&e^{\lambda_{1}t}\langle |\alpha \widetilde{U}^{\alpha}\phi-\phi|h,X_{t}\rangle+\left|e^{\lambda_{1}t}\langle (\alpha \widetilde{U}^{\alpha} \phi) h,X_{t}\rangle-(\phi h,\hat{h}) W^{h}_{\infty}(X)\right|\nonumber\\
&\le&\|\alpha \widetilde{U}^{\alpha}\phi-\phi\|_{\infty}W^{h}_{t}(X)+\left|e^{\lambda_{1}t}\langle (\alpha \widetilde{U}^{\alpha} \phi) h,X_{t}\rangle-(\phi h,\hat{h}) W^{h}_{\infty}(X)\right|.\nonumber
\end{eqnarray}
By Lemmas \ref{lem7.5} and \eqref{them7.3.2}, we conclude \eqref{them7.3.1} by letting
first $t\to{\infty}$ and then $\alpha\to{\infty}$.\qed

\medskip
{\bf Sandra Palau}

Department of Statistics and Probability, Instituto de Investigaciones en Matem\'aticas Aplicadas y en 

Sistemas, Universidad Nacional Aut\'onoma de M\'exico, M\'exico.

Email address: sandra@sigma.iimas.unam.mx

\medskip
{\bf Ting Yang}

School of Mathematics and Statistics, Beijing Institute of Technology, Beijing, 100081, P.R.China;

Beijing Key Laboratory on MCAACI, Beijing, 100081, P.R. China.

Email address: yangt@bit.edu.cn

\end{document}